\documentclass{article}

\PassOptionsToPackage{numbers, compress}{natbib}

\usepackage[utf8]{inputenc} 
\usepackage[T1]{fontenc}    
\usepackage{hyperref}       
\usepackage{url}            
\usepackage{booktabs}       
\usepackage{amsfonts}       
\usepackage{nicefrac}       
\usepackage{microtype}      
\usepackage{xcolor}         
\usepackage{mleftright}

\RequirePackage{graphicx}
\RequirePackage{amsthm,amsmath,amsfonts,amssymb}
\RequirePackage[numbers]{natbib}

\newcommand{\bm}{\boldsymbol}

\usepackage{verbatim}
\newtheorem{theorem}{Theorem}[section]
\newtheorem{lemma}[theorem]{Lemma}
\newtheorem{corollary}[theorem]{Corollary}
\newtheorem{proposition}[theorem]{Proposition}
\newtheorem{assumption}{Assumption}

\usepackage[normalem]{ulem}

\usepackage{PRIMEarxiv}
\usepackage{footmisc} 
\usepackage[utf8]{inputenc} 
\usepackage[T1]{fontenc}    
\usepackage{hyperref}       
\usepackage{url}            
\usepackage{booktabs}       
\usepackage{amsfonts}       
\usepackage{nicefrac}       
\usepackage{microtype}      
\usepackage{lipsum}
\usepackage{fancyhdr}       
\usepackage{graphicx}       
\graphicspath{{media/}}     
\usepackage{hyperref}
\usepackage{url}
\usepackage{graphicx}
\usepackage{subfigure}
\usepackage{mathrsfs}

\usepackage{comment}

\fancypagestyle{firstpage}
{
    \lfoot{12th International Conference on Learning Representations (ICLR 2024)}
    \cfoot{\quad}
}

\title{The Optimality of Kernel Classifiers in Sobolev Space}

\author{Jianfa Lai\footnotemark[1]\\
Tsinghua University, Beijing, China \\
\texttt{Jianfalai@mail.tsinghua.edu.cn, } \\
\And
 Zhifan Li\footnotemark[1]\\
Department of Mathematical Sciences  \\
\texttt{stahd@nus.edu.sg}\\
\And
Dongming Huang\\
National University of Singapore, Singapore \\
\texttt{stahd@nus.edu.sg}\\
\And
Qian Lin\footnotemark[2] \footnotemark[3]\\
Tsinghua University, Beijing, China \\
\texttt{qianlin@tsinghua.edu.cn}
}

\begin{document}
\maketitle

\newtheorem{example}{Example}
\newtheorem{definition}{Definition}
\newtheorem{conjecture}{Conjecture}
\newtheorem{remark}{Remark}
\renewcommand{\thefootnote}{\fnsymbol{footnote}} 
\footnotetext[1]{These authors contributed equally to this work.} 
\footnotetext[3]{Qian Lin also affiliates with Beijing Academy of Artificial Intelligence, Beijing, China.} 

\begin{abstract}
Kernel methods are widely used in machine learning, especially for classification problems. However, the theoretical analysis of kernel classification is still limited. This paper investigates the statistical performances of kernel classifiers. With some mild assumptions on the conditional probability $\eta(x)=\mathbb{P}(Y=1\mid X=x)$, we derive an upper bound on the classification excess risk of a kernel classifier using recent advances in the theory of kernel regression. We also obtain a minimax lower bound for Sobolev spaces, which shows the optimality of the proposed classifier. Our theoretical results can be extended to the generalization error of overparameterized neural network classifiers. 
To make our theoretical results more applicable in realistic settings, we also propose a simple method to estimate the interpolation smoothness of $2\eta(x)-1$ and apply the method to real datasets.
\end{abstract}

\thispagestyle{firstpage}

\section{Introduction}

In this paper, we study the problem of binary classification in a reproducing kernel Hilbert space (RKHS).  
Suppose $n$ i.i.d samples $\{(X_i,Y_i)\in \mathcal{X}\times \{-1,1\}\}$  are drawn from a joint distribution $(X,Y)\sim \rho$, where the conditional probability of the response variable $Y$ given the predictor variable $X=x$ is denoted by $\eta(x)= \mathbb{P}(Y=1|X=x)$. 
We aim to find a classifier function $f(x) : \mathcal{X} \to [-1,1]$ that minimizes the classification risk, defined as:
$$
\mathcal{L}(\hat{f}):= \mathbb{P}_{(X,Y)\sim \rho)}\left[\operatorname{sign}(\hat{f}(X)) \neq Y\right] .
$$
The minimal classification risk is achieved by the Bayes classifier function corresponding to $\rho$, which is defined as $f^*_{\rho}(x) = 2\eta(x)-1$. 
Our main focus is on analyzing the convergence rate of the classification excess risk, defined as:
$$
\mathcal{E}(\hat{f}) = \mathcal{L}(\hat{f}) - \mathcal{L}(f^*_{\rho}).
$$
 
This paper studies a class of kernel methods called spectral algorithms (which will be defined in Section~\ref{spectra}) for constructing estimators of $f_\rho^*$. The candidate functions are selected from an RKHS $\mathcal{H}$, which is a separable Hilbert space associated with a kernel function $K$ defined on $\mathcal{X}$ \citep{smale2007learning,steinwart2008support}. 
Spectral algorithms, as well as kernel methods, are becoming increasingly important in machine learning because both experimental and theoretical results show that overparameterized neural network classifiers exhibit similar behavior to classifiers based on kernel methods \cite{belkin2018understand}. Therefore, understanding the properties of classification with spectral algorithms can shed light on the generalization of deep learning classifiers.


In kernel methods context, researchers often assume that $f_\rho^* \in \mathcal{H} $, and have obtained the minimax optimality of spectral algorithms \citep{caponnetto2006optimal,caponnetto2007optimal}. Some researchers have also studied the convergence rate of the generalization error of misspecified spectral algorithms $\left(f_\rho^* \notin \mathcal{H} \right)$, assuming that $f_\rho^*$ falls into the interpolation space $[ \mathcal{H} ]^s$ with some $s>0$ \citep{fischer2020sobolev,zhang2023optimality}. 
In this line of work, researchers consider the embedding index condition which reflects the capability of $\mathcal{H}$ embedding into $L^{\infty}$ space. Moreover, 
\cite{zhang2023optimality} extends the boundedness assumption to the cases where $[ \mathcal{H} ]^s \cap L^{\infty}( X , \mu) \varsubsetneqq[ \mathcal{H} ]^s$.

Motivated by the aforementioned studies, we adopt similar assumptions in our study of kernel classifiers trained via the gradient flow. 
We assume that the Bayes classifier $f_\rho^* \in [\mathcal{H}]^s$  satisfies the boundedness condition $f^*_{\rho} \in [-1,1]$. 
We first derive the upper bound of the classification excess risk, showing that the generalization error of the kernel classifier is highly related to the interpolation smoothness $s$. 
To clarify the minimax optimality of kernel classification, we then obtain the minimax lower bound for classification in Sobolev RKHS, which is a novel result in the literature. 
Our technique is motivated by the connection between kernel estimation and infinite-width neural networks, and our framework can be applied to neural network supervised learning. Furthermore, we provide a method to estimate the interpolation space smoothness parameter $s$ and also present some numerical results for neural network classification problems through simulation studies and real data analysis.

\subsection{Our contribution}

In this paper, we study the generalization error of kernel classifiers. We show that
\begin{itemize}
    \item [$i)$] We show the generalization error of the gradient flow kernel classifier is bounded by $O(n^{- s\beta/(2 s\beta+2)})$ provided that the Bayes classifier $f_\rho^*\in [\mathcal{H}]^{s}$, where $\beta$ is the eigenvalue decay rate (EDR) of the kernel. 
    This result is not only applicable to the Sobolev RKHS $\mathcal{H}$ but also to any RKHS with the embedding index $\alpha_0=1/\beta$, such as the RKHS with dot-product kernels and the RKHS with shift-invariant periodic kernels.
    \item [$ii)$] We establish a minimax lower bound on the classification excess risk in the interpolation space of Sobolev RKHS. Combined with the results in $i)$, the convergence rate of the kernel classifier is minimax optimal in Sobolev space. Before our work, \cite{yang1999minimax} illustrated a similar result of the minimax lower bound for Besov spaces. However, the result has only been proved for $d=1$  by \cite{kerkyacharian1992density} and the case for $d>1$ remains unresolved. 
    \item [$iii)$] To make our theoretical results more applicable in realistic settings, we propose a simple method to estimate the interpolation smoothness $s$. 
    We apply this method to estimate the relative smoothness of various real datasets with respect to the neural tangent kernels, where the results are in line with our understanding of these real datasets.
\end{itemize}




\subsection{Related works}

We study the classification rules derived from a class of real-valued functions in a reproducing kernel Hilbert space (RKHS), which are used in kernel methods such as Support Vector Machines (SVM) \citep{steinwart2008support}. Most of the existing works consider hinge loss as the loss function, i.e. \cite{wahba2002soft,steinwart2007fast, bartlett2008classification,blanchard2008statistical} etc. Another kernel method, kernel ridge regression, also known as least-square SVM  \citep{steinwart2008support}, is investigated by some researchers \citep{xiang2009classification,rifkin2003regularized}. Recently, some works have combined the least square loss classification with neural networks \citep{demirkaya2020exploring,hu2021understanding}.

We choose kernel methods because it allows us to use the integral operator tool for analysis \citep{de2005model,caponnetto2007optimal,fischer2020sobolev,zhang2023optimality}, while previous SVM works tend to use the empirical process technique \citep{steinwart2007fast}. Moreover, we can easily extend the $\, \mathcal{H}$ to the misspecified model case $[\mathcal{H}]^s$ when true model $f^*_{\rho}$ belongs to a less-smooth interpolation space. Furthermore, we consider more regularization methods, collectively known as spectral algorithms, which were first proposed and studied by \cite{rosasco2005spectral,bauer2007regularization,caponnetto2007optimal}. \cite{zhang2023optimality} combined these two ideas and obtained minimax optimality for the regression model. We extend their results to the classification problems.

We study the minimax optimality of Sobolev kernel classification, and before our work, the minimax lower bound of classification excess risk for the RKHS class was seldom considered. \cite{loustau2008aggregation,loustau2009penalized} have discussed Classification problems in Sobolev space, but they did not consider the lower bound of classification risk. \cite{audibert2004classification,audibert2007fast,massart2006risk} provided some minimax lower bound techniques for classification, but how to solve RKHS remains
unknown. Sobolev space (see, e.g., \cite{adams2003sobolev}) is known as a vector space of functions equipped with a norm that is a combination of $L^2$-norms of the function together with its derivatives up to a given order and can be embedded into H\"{o}lder class. Inspired by the minimax lower bound for H\"{o}lder class classification in \cite{audibert2007fast}, we derive the lower bound for the Sobolev class.  

Recently, deep neural networks have gained incredible success in classification tasks from image classification \citep{krizhevsky2012imagenet,he2016deep} to natural language processing \citep{devlin2019bert}. Since \cite{jacot2018neural} introduced the neural tangent kernel,  The gradient flow of the training process can be well approximated by a simpler gradient flow associated with the
NTK kernel when the width of neural networks is sufficiently large \citep{lai2023generalization,li2023statistical}. 
Therefore, we can analyze the classification risk of neural networks trained by gradient descent.

\section{Preliminaries}

We observe $n$ samples $\{(X_i,Y_i)\in \mathcal{X}\times \{-1,1\}\}$ where $\mathcal{X}\subset \mathbb{R}^{d}$ is compact. Let $ \rho $ be an unknown probability distribution on $\mathcal{X} \times \{-1,1\}$ and $\mu$ be the marginal distribution on $\mathcal{X}$. We assume $\mu$ has a uniformly bounded density $0<\mu_{min} \leq\mu(x)\leq \mu_{max}$ for $x\in \mathcal{X}$. The classification task is to predict the unobserved label $y$ given a new input $x$. The conditional probability is defined as  $\eta(x)= \mathbb{P}(Y=1|X=x)$. For any classifier $f$, the risk based on the 0-1 loss can be written as 
\begin{equation}
    \mathcal{L}(f)=\bm{E}_{(X,Y)\sim \rho}\mathbb{I}\{\operatorname{sign} (f(X))\neq Y\} = \bm{E}_{X} [(1-\eta(X))\mathbb{I}\{f(X)\geq 0\} + \eta(X) \mathbb{I}\{f(X)<0\}].
\end{equation}
One of the minimizers of the risk has the form $f^*_\rho=2\eta-1$. Let $\mathcal{L}^*=\mathcal{L}(f^*_{\rho})$. 
For any classifier $\hat{f}$ learned from data, its accuracy is often characterized by the classification excess risk, which can be formulated as 
\begin{equation}\label{excessrisk}
    \begin{aligned}
    &\mathcal{E}(\hat{f}) = \mathcal{L}(\hat{f}) - \mathcal{L}^* =\bm{E}_{X}(|f^*_{\rho}(X)| \mathbb{I}\{\hat{f}(X)f^*_{\rho}(X)<0\} ).
\end{aligned}
\end{equation}
In the rest of this section, we introduce some essential concepts in RKHS and kernel classifiers. In Section 2.1, we review some definitions in the interpolation space of RKHS. The relationship between fractional Sobolev space and Sobolev RKHS is presented in Section 2.2. Section 2.3 presents the explicit formula of the gradient-flow kernel classifier and the corresponding rewritten form through spectral algorithms and filter functions.

\subsection{Interpolation Space of RKHS}

Denote $L^2(\mathcal{X}):=L^2(\mathcal{X},\mu)$ as the $L^2$ space. Throughout the paper, we denote by $\mathcal{H}$ a separable RKHS on $\mathcal{X}$ with respect to a continuous kernel function $K$. We also assume that $\sup_{x\in\mathcal{X}}K(x,x)\leq \kappa$ for some constant $\kappa$.
The celebrated Mercer's theorem shows that there exist non-negative numbers $\lambda_{1} \geq \lambda_{2} \geq \cdots$ and functions $e_{1}, e_{2}, \cdots \in L^{2}(\mathcal{X})$ such that $\left<e_{i},e_{j}\right>_{L^{2}(\mathcal{X})}=\delta_{ij}$ and 
 \begin{align}
K_{d}(x,x')=\sum_{j=1}^{\infty}\lambda_{j}e_{j}(x)e_{j}(x'),
 \end{align}
 where the series on the right hand side converges in $L^{2}(\mathcal{X})$. 
 
Denote the natural embedding inclusion operator by $S_k: \mathcal{H} \rightarrow L^2( \mathcal{X} , \mu)$. Moreover, the adjoint operator $S_k^*: L^2( \mathcal{X} , \mu) \rightarrow \mathcal{H} $ is an integral operator, i.e., for $f \in L^2( \mathcal{X} , \mu)$ and $x \in \mathcal{X}$, we have
$$
\left(S_k^* f\right)(x)=\int_{ \mathcal{X} } K\left(x, x^{\prime}\right) f\left(x^{\prime}\right) d \mu\left(x^{\prime}\right) .
$$
It is well-known that $S_k$ and $S_k^*$ are Hilbert-Schmidt operators (and thus compact) and their HS norms (denoted as $\|\cdot\|_2$ ) satisfy that
$$
\left\|S_k^*\right\|_2=\left\|S_k\right\|_2=\|K\|_{L^2( \mathcal{X} , \mu)}:=\left(\int_{ \mathcal{X} } K(x, x) d \mu(x)\right)^{1 / 2} \leq \kappa .
$$
Next, we define two integral operators as follows:
$$
L:=S_k S_k^*: L^2( \mathcal{X} , \mu) \rightarrow L^2( \mathcal{X} , \mu), \quad T:=S_k^* S_k: \mathcal{H}  \rightarrow \mathcal{H}  .
$$
$L$ and $T$ are self-adjoint, positive-definite, and in the trace class (and thus Hilbert-Schmidt and compact). Their trace norms (denoted as $\|\cdot\|_1$ ) satisfy that $\left\|L\right\|_1=\|T\|_1=\left\|S_k\right\|_2^2=\left\|S_k^*\right\|_2^2$.


For any $ s \ge 0$, the fractional power integral operator $L^{s}: L^{2}(\mathcal{X},\mu) \to L^{2}(\mathcal{X},\mu)$ and $T^{s}: \mathcal{H} \to \mathcal{H} $ are defined as 
\begin{align}
  L^{s}(f)=\sum_{j=1}^{\infty} \lambda_j^{s} \left\langle f, e_j\right\rangle_{L^2} e_j, \quad T^{s}(f)=\sum_{j=1}^{\infty} \lambda_j^{s} \left\langle f, \lambda_{j}^{\frac{1}{2}}e_j\right\rangle_{\mathcal{H}}\lambda_{j}^{\frac{1}{2}} e_j.
\end{align}
The interpolation space $[\mathcal{H}]^s $ is defined as
\begin{align}
  [\mathcal{H}]^s :=  \left\{\sum_{j =1}^{\infty} a_j \lambda_j^{s / 2}e_{j}: \sum_{j}^{\infty} a_j^2<\infty \right\} \subseteq L^{2}(\mathcal{X})
\end{align}

It is easy to show that $[\mathcal{H}]^s $ is also a separable Hilbert space with orthogonal basis $ \{ \lambda_{i}^{s/2} e_{i}\}_{i \in N}$. Specially, we have $[\mathcal{H}]^0 \subseteq L^{2}(\mathcal{X},\mu) $, $[\mathcal{H}]^1 = \mathcal{H}$ and $[\mathcal{H}]^{s_{2}} \subsetneq [\mathcal{H}]^{s_{1}} \subsetneq [\mathcal{H}]^0 $ for any numbers $0 < s_{1} < s_{2}$. 
For the functions in $[\mathcal{H}]^{s}$ with larger $s$, we say they have higher (relative) interpolation smoothness with respect to the RKHS (the kernel).

\subsection{Fractional Sobolev Space and Sobolev RKHS}
 For $m\in\mathbb{N}$, we denote the usual Sobolev space $W^{m,2}(\mathcal{X})$ by $H^{m}(\mathcal{X})$ and $L^2(\mathcal{X})$ by $H^{0}(\mathcal{X})$. Then the (fractional) Sobolev
space for any real number $r>0$ can be defined through the real interpolation 
\begin{align*}
    H^{r}(\mathcal{X}):= \left( L^2(\mathcal{X}), H^{m}(\mathcal{X}) \right)_{\frac{r}{m},2}
\end{align*}
where $m:=min\{k\in\mathbb{N} \ :\ k>r \}$. 

It is well known that when $r>d/2$, $H^{r}$ is a separable RKHS with respect to a bounded kernel and the corresponding eigenvalue decay rate (EDR) is $\beta=2r/d$. Furthermore, the interpolation space of $H^{r}(\mathcal{X})$ under Lebesgue measure is given by
\begin{equation}\label{sobolev_interpolation}
  [H^{r}(\mathcal{X})]^{s} = H^{rs}(\mathcal{X}).
\end{equation}

It follows that given a Sobolev RKHS $\mathcal{H}=H^{r}$ for $r>d/2$, if  $f\in H^{a}$ for any $a>0$, one can find that $f\in [\mathcal{H}]^s$ with $s=a/r$. Thus, in this paper, we will assume that the Bayes classifier $f_{\rho}^{*}$
is in the interpolation of the Sobolev RKHS $[\mathcal{H}]^s$.

\subsection{Kernel Classifiers: Spectra Algorithm}\label{spectra}


%
 We then introduce a more general framework known as spectra algorithm \citep{rosasco2005spectral,caponnetto2006optimal,bauer2007regularization}. We define the filter function and the spectral algorithms as follows:

\begin{definition}[Filter function] 
Let $\left\{\varphi_\nu:\left[0, \kappa^2\right] \rightarrow R ^{+} \mid \nu \in \Gamma \subseteq R ^{+}\right\}$ be a class of functions and $\psi_\nu(z)=1-z \varphi_\nu(z)$. If $\varphi_\nu$ and $\psi_\nu$ satisfy:
\begin{itemize}
    \item  $\forall \alpha \in[0,1]$, we have
 $\quad\label{filter1}
    \sup _{z \in\left[0, \kappa^2\right]} z^\alpha \varphi_\nu(z) \leq E \nu^{1-\alpha}, \quad \forall \nu \in \Gamma ;
$

\item $\exists \tau \geq 1$ s.t. $\forall \alpha \in[0, \tau]$, we have
$\quad\label{filter2}
    \sup _{z \in\left[0, \kappa^2\right]}\left|\psi_\nu(z)\right| z^\alpha \leq F_\tau \nu^{-\alpha}, \quad \forall \nu \in \Gamma,
$
\end{itemize}

where $E, F_\tau$ are absolute constants, then we call $\varphi_\nu$ a filter function. We refer to $\nu$ as the regularization parameter and $\tau$ as the qualification.
\end{definition}

\begin{definition}[spectral algorithm] Let $\varphi_\nu$ be a filter function index with $\nu>0$. Given the samples $Z$, a spectral algorithm produces an estimator of $f_\rho^*$ given by
$
\hat{f}_\nu=\varphi_\nu\left(T_X\right) g_Z.
$
\end{definition}

The following example shows that $\hat{f}_t(x)$ can be formulated by the spectral algorithms.
\begin{example}[Classifier with Gradient flow]\label{eq:gradient_flow}
The filter function of gradient flow $\varphi_\nu$ can be defined as
$
\varphi_\nu^{ gf }(z)=\frac{1-e^{-\nu z}}{z}.
$
The qualification $\tau$ could be any positive number, $E=1$ and $F_\tau=(\tau / e)^\tau$. So that for a test input $x$, the predicted output is given by $\hat{y} = \operatorname{sign} (\hat{f}_{\nu}(x))$. 
\end{example} 

 Other spectral algorithms consist of kernel ridge regression, spectral cut-off, iterated Tikhonov, and so on. For more examples, we refer to \cite{gerfo2008spectral}. Spectral algorithms differ in $\varphi_\nu(z)$ and $\psi_{\nu}(z)$, which is corresponding to saturation effect defined in \cite{gerfo2008spectral}. Moreover, \cite{li2023_SaturationEffect} gives a thorough analysis of the saturation effect for kernel ridge regression.

\paragraph{Notations.}
Denote $B(x,r)$ as a ball, and $\lambda[B(x,r)]$ is denoted as the Lebesgue measure of $B(x,r)$.  We use $\|\cdot\|_{ \mathscr{B}  \left(B_1, B_2\right)}$ to denote the operator norm of a bounded linear operator from a Banach space $B_1$ to $B_2$, i.e., $\|A\|_{ \mathscr{B} \left(B_1, B_2\right)}=\sup _{\|f\|_{B_1=1}}\|A f\|_{B_2}$. Without bringing ambiguity, we will briefly denote the operator norm as $\|\cdot\|$. In addition, we use $\operatorname{tr} A$ and $\|A\|_1$ to denote the trace and the trace norm of an operator. We use $\|A\|_2$ to denote the Hilbert-Schmidt norm.

\section{Main Results}


\subsection{Assumptions}
This subsection lists the standard assumptions for general RKHS $\mathcal{H}$ and clarifies how these assumptions correspond to properties of Sobolev RKHS.

\begin{assumption}[Source condition]\label{ass:rkhs}
  For $s > 0 $, there is a constant $B > 0 $ such that $f_{\rho}^{*} \in [\mathcal{H}]^{s}$ and $\| f_{\rho}^{*} \|_{[\mathcal{H}]^{s}} \le B$.
 
\end{assumption}
This assumption is weak since $s$ can be small. However, functions in $[\mathcal{H}]^{s}$ with smaller $s$ are less smooth, which will be harder for an algorithm to estimate.
\begin{assumption}[Eigenvalue Decay Rate (EDR)]\label{ass:eigendecay_rate}
    The EDR of the eigenvalues $\{\lambda_j\}$ associated to the kernel $K$ is $\beta>1$, i.e.,
    \begin{align}
        cj^{-\beta} \leq \lambda_j \leq Cj^{-\beta}
    \end{align}
    for some positive constants $c$ and $C$. 
\end{assumption}
Note that the eigenvalues $\lambda_i$ and EDR are only determined by the marginal distribution $\mu$ and the RKHS $\mathcal{H}$. For Sobolev RKHS $H^r$ equipped with  Lebesgue measure $\nu$ and bounded domain with smooth boundary $\mathcal{X} \subseteq R ^d$, it is well known that when $r>d/2$, $H^{r}$ is a separable RKHS with respect to a bounded kernel and the corresponding eigenvalue decay rate (EDR) is $\beta=2r/d$ \citep{edmunds_triebel_1996}.

Our next assumption is the embedding index. First, we give the definition of embedding property \cite{fischer2020sobolev}: For $0<\alpha < 1$, there is a constant $A>0$ with
$
\left\|[H]_\nu^\alpha \hookrightarrow L_{\infty}(\nu)\right\| \leq A .
$
This means $[\mathcal{H}]^\alpha$ is continuously embedded into $L_{\infty}(\nu)$ and the operator norm of the embedding is bounded by $A$. The larger $\alpha$ is, the weaker the embedding property is.

\begin{assumption}[Embedding index]\label{ass:emb}
Suppose that there exists $\alpha_{0} > 0$, such that
\begin{displaymath}
  \alpha_{0} = \inf\left\{ \alpha \in [\frac{1}{\beta},1] :  \left\|[\mathcal{H}]^\alpha \hookrightarrow L^{\infty}(\mathcal{X},\mu)\right\| < \infty  \right\},
\end{displaymath}
and we refer to $\alpha_{0}$ as the \textit{embedding index} of an RKHS $\mathcal{H}$. 
\end{assumption}

This assumption directly implies that all the functions in $[\mathcal{H}]^{\alpha}$ are $\mu$-a.e bounded for $\alpha>\alpha_0$.  Moreover, we will clarify this assumption for Sobolev kernels and dot-product kernels on $\mathbb{S}^{d-1}$ in the appendix.

\subsection{Minimax optimality of kernel classifiers}
This subsection presents our main results on the minimax optimality of kernel classifiers. We first establish a minimax lower bound for the Sobolev RKHS $H^r(\mathcal{X})$ under the source condition (Assumption \ref{ass:rkhs}). We then provide an upper bound based on Assumptions \ref{ass:rkhs}, \ref{ass:eigendecay_rate}, and \ref{ass:emb}, and we clarify that the Sobolev RKHS satisfies these assumptions. As a result, we demonstrate that the Sobolev kernel classifier is minimax rate optimal.

\begin{theorem}[Lower Bound]\label{thm:lower_bound}
 Suppose $f^*_\rho\in [H^r(\mathcal{X})]^s$ for $s>0$, where $H^r$ is the Sobolev RKHS. 
 For all learning methods $\hat{f}$, for any fixed $\delta \in (0,1)$, when $n$ is sufficiently large, there is a distribution $\rho \in \mathcal{P}$ such that, with probability at least $1 - \delta$, we have 
   \begin{equation}
       \mathcal{E}(\hat{f})  \ge C \delta n^{-\frac{ s\beta}{ 2(s\beta+1)}}, 
   \end{equation}
   where $C$ is a universal constant. 
\end{theorem}

Theorem \ref{thm:lower_bound} shows the minimax lower bound on the classification excess risk over the interpolation space of the Sobolev RKHS. 
Theorem \ref{thm:lower_bound} also establishes a minimax lower bound at the rate of $n^{-\frac{a}{2a+d}}$ for the Sobolev space $H^{a}$ with $a=rs$.  
\cite{yang1999minimax} illustrated a similar result of the minimax lower bound for Besov spaces. However, the result has only been proved for $d=1$  by \cite{kerkyacharian1992density} and the case for $d>1$ remains unresolved. 

The following theorem presents an upper bound for the kernel classifier.
\begin{theorem}[Upper Bound]\label{thm:upper bound}
    Suppose that Assumptions  \ref{ass:rkhs}, \ref{ass:eigendecay_rate}, and \ref{ass:emb} hold for $0<s\leq 2\tau$, where $\tau$ is the qualification of the filter function. By choosing $\nu \asymp n^{\frac{ \beta }{  s\beta + 1}}$, for any fixed $\delta \in (0,1)$, when $n$ is sufficiently large, with probability at least $1 - \delta$, we have
    \begin{align}
        \mathcal{E}(\hat{f}_\nu) \leq C \left(\ln \frac{4}{\delta}\right) n^{-\frac{ s\beta}{ 2(s\beta+1)}}
    \end{align}
     where $C$ is a constant independent of $n$ and $\delta$.
\end{theorem}

%



Combined with Theorem \ref{thm:lower_bound}, Theorem \ref{thm:upper bound} shows that by choosing a proper early-stopping time, the Sobolev kernel classifier is minimax rate optimal. Moreover, given the kernel and the decay rate $\beta$, the optimal rate is mainly affected by the smoothness $s$ of $f_\rho^*$ with respect to the kernel. Thus, in Section \ref{sec:smoothness}, we will introduce how to estimate the smoothness of functions or datasets given a specific kernel.

We emphasize that Theorem \ref{thm:upper bound} can be applied to any general RKHS with an embedding index $\alpha_0=1/\beta$, such as an RKHS with a shift-invariant periodic kernel and an RKHS with a dot-product kernel. Thanks to the uniform convergence of overparameterized neural networks \citep{lai2023generalization,li2023statistical}, Theorem \ref{thm:upper bound} can also be applied to analyze the generalization error of the neural network classifiers. We will discuss this application in the next section.

\section{Applications in Neural Networks}
Suppose that we have observed $n$ i.i.d. samples $\{X_i,Y_i\}_{i=1}^{n}$ from $\rho$. For simplicity, we further assume that the marginal distribution $\mu$ of $\rho$ is the uniform distribution on the unit sphere $\mathbb{S}^{d-1}$. We use a neural network with $L$ hidden layers and width $m$ to perform the classification on $\{X_i,Y_i\}_{i=1}^{n}$. The network model $ f(x;\theta)$ and the resulting prediction are given by the following equations
\begin{align}
  \begin{split}
    h^{0}(x) & = x,\quad   h^{l}(x)= \sqrt{\frac{2}{m}}\sigma(W^{l-1}h^{l-1}(x)),\quad l=1,...,L\\
    f(x;\theta) & = W^{L}h^{L}(x)\quad \text{and} \quad \hat{y} = \operatorname{sign}(f(x;\theta)),
  \end{split}
\end{align}
where $h^{l}$ represents the hidden layer, $\sigma(x) := \max(x,0)$ is the ReLU activation (applied elementwise), $W^{0}\in \mathbb{R}^{m\times d}$ and $W^l\in \mathbb{R}^{m\times m}$ are the parameters of the model.
We use $\theta$ to represent the collection of all parameters flatten as a column vector. With the mirrored initialization (shown in \cite{li2023statistical}), we consider the training process given by the gradient flow $\dot{\theta} = -\partial L(\theta)/ \partial \theta$, where the squared loss function is 
adopted $L(\theta) = \frac{1}{2n} \sum_{i=1}^{n} \left(Y_i - f(X_i,\theta)) \right)^2. $

The consideration for this choice of loss function is that the squared loss function is robust for optimization and more suitable for hard learning scenarios (\cite{hui2020evaluation,demirkaya2020exploring,kornblith2020demystifying}). \cite{hui2020evaluation} showed that the square loss function has been shown to perform well in modern classification tasks, especially in natural language processing while \cite{kornblith2020demystifying} presented the out-of-distribution robustness of the square loss function.


When the network is overparameterized, \cite{li2023statistical} showed that the trained network $f(x;\theta)$ can be approximated by a kernel gradient method with respect to the following neural tangent kernel
\begin{align}
  \label{eq:NTK_Formula}
  K_{ntk}(x,x') =
   \sum_{r=0}^L \kappa^{(r)}_1(\bar{u}) \prod_{s=r}^{L-1} \kappa_0(\kappa^{(s)}_1(\bar{u})) 
\end{align}
where $\bar{u} = \langle x,x'\rangle$,
$\kappa_1^{(p)} = \kappa_1\underbrace{\circ \dots \circ}_{p \text{ times}} \kappa_1$ represents $p$ times composition of $\kappa_1$ and $\kappa_1^{(0)}(u) = u$ by convention; if $r = L$, the product $\prod_{s=r}^{L-1}$ is understood to be $1$. Denote $ Y_{[n]} = (Y_1,...,Y_n)^{T}$, $K(X_{[n]},X_{[n]})$ as an $n\times n$ matrix of $\left(K(X_{i},X_{j})  \right)_{i,j\in[n]}$ and $\lambda_{min} = \lambda_{min}(K(X_{[n]},X_{[n]}))$
The following proposition shows the uniform convergence of $f(x;\theta)$.

\begin{proposition}[Theorem 1 in \cite{li2023statistical}]
  \label{prop:UnifConverge}
  Suppose $x\in\mathbb{S}^{d-1}$. For any $\epsilon > 0$, 
  any hidden layer $L\geq 2$,
  and $\delta \in (0,1)$,
  when the width $m \geq \operatorname{poly}\left(n, \lambda_{min}^{-1}, || Y_{[n]} ||_2, \ln (1 / \delta), \ln(1/\epsilon) \right)$, with probability at least $1-\delta$ with respect to random initialization, we have
  \begin{align*}
    \sup_{t\geq 0} \sup_{x \in \mathcal{X}} |f_t(x;\theta) - f_t^{ntk}(x)| \leq \epsilon.
  \end{align*}
  where $f_{t}^{ntk}(x)$ is defined as in Example \ref{eq:gradient_flow} but with the kernel $K_{ntk}$.
\end{proposition}

Theorem G.5 in \cite{haas2023mind} showed that the RKHS of the NTK on $\mathbb{S}^{d-1}$ is a Sobolev space. Moreover, the kernel $K_{ntk}$ is a dot-product kernel satisfying a polynomial eigenvalue decay $\beta = d/(d-1)$. Thus, we can obtain the following corollary by combining Theorem \ref{thm:upper bound} and Proposition \ref{prop:UnifConverge}. 
\begin{corollary}\label{coro:upper bound}
   Suppose that $x\in\mathbb{S}^{d-1}$ and Assumption  \ref{ass:rkhs} holds for $\mathcal{H}$ being the RKHS of the kernel $K_{ntk}$  and $s>0$. 
    Suppose $t \asymp n^{\frac{ \beta }{  s\beta + 1}}$. For any fixed $\delta \in (0,1)$, when $m \geq \operatorname{poly}\left(n, \lambda_{min}^{-1}, || Y_{[n]} ||_2, \ln (1 / \delta) \right)$ and $n$ is sufficiently large, with probability at least $1 - \delta$, we have
    \begin{align}
        \mathcal{E}(f_t(x;\theta)) \leq C \left(\ln \frac{4}{\delta}\right) n^{-\frac{ s\beta}{ 2(s\beta+1)}}
    \end{align}
     where $C$ is a constant independent of $n$ and $\delta$.
\end{corollary}

This corollary shows that the generalization error of a fine-tuned, overparameterized neural network classifier converges at the rate of $n^{-\frac{ s\beta}{ 2(s\beta+1)}}$. 
This result also highlights the need for additional efforts to understand the smoothness of real datasets with respect to the neural tangent kernel. 
A larger value of $s$ corresponds to a faster convergence rate, indicating the possibility of better generalization performance. 
Determination of the smoothness parameter $s$ will allow us to assess the performance of an overparameterized neural network classifier on a specific dataset.

\section{Estimation of smoothness}\label{sec:smoothness}
In this section, we provide a simple example to illustrate how to determine the relative smoothness $s$ of the ground-truth function with respect to the kernel. 
Then we introduce a simple method to estimate $s$ with noise and apply the method to real datasets with respect to the NTK. 

\paragraph{Determination of $s$.} Suppose that $\mathcal{X}\in[0,1]$ and the marginal distribution $\mu_{\mathcal{X}}$ is a uniform distribution on $[0,1]$. We consider the min kernel $K_{min}(x,x')=\min(x,x')$ \citep{wainwright2019high} and denote by $\mathcal{H}_{min}$ the corresponding RKHS. 
The eigenvalues and the eigenfunctions of $\mathcal{H}_{min}$ are 
\begin{align}
    \lambda_j=\left(\frac{2j-1}{2}\pi\right)^{-2}, \quad e_j(x) =\sqrt{2}\sin(\frac{2j-1}{2}\pi x), \quad j\geq 1.
\end{align}
Thus, the EDR is $\beta=2$. 
For illustration, we consider the ground true function $ f^*(x) = \cos(2\pi x)$. 
Suppose $f^*(x) = \sum_j^{\infty} f_j e_j(x)$, then we have $f_j = \sqrt{2} \int_{0}^{1} \cos(2\pi x) \sin(\frac{2j-1}{2}\pi x)dx \asymp j^{-1}$. Thus, $f_j\asymp j^{-r}$ where $r=1$. By the definition of the interpolation space, we have  $ s=\frac{2r-1}{\beta} = 0.5$.

\paragraph{Estimation of $s$ in regression. }
\begin{figure}[htb]
  \centering
\subfigure[ ]{
    \begin{minipage}[t]{0.3\linewidth}
      \centering
      \includegraphics[width=1\linewidth]{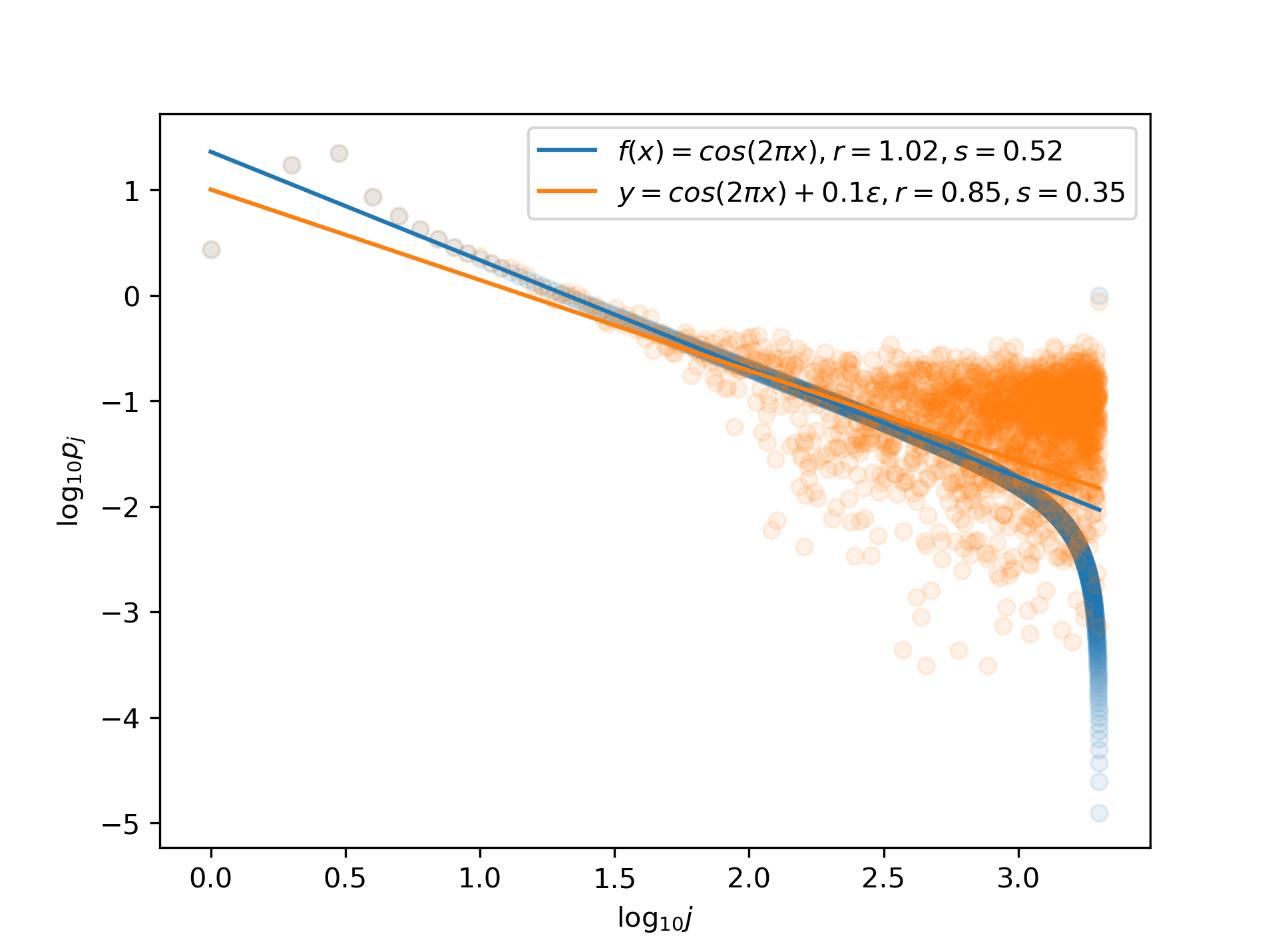}
    \end{minipage}
    }
\subfigure[ ]{
    \begin{minipage}[t]{0.3\linewidth}
      \centering
      \includegraphics[width=1\linewidth]{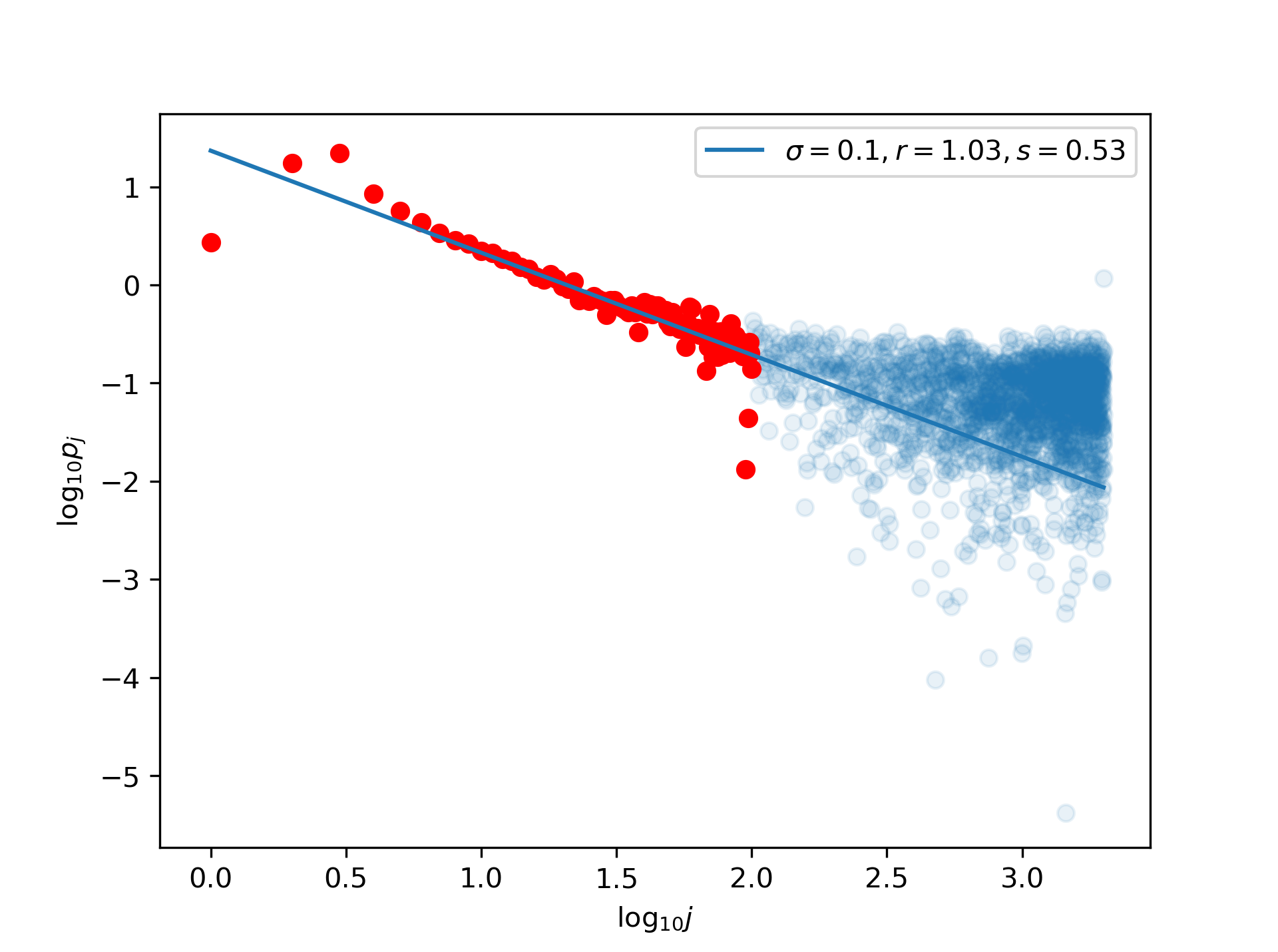}
    \end{minipage}%
}
\subfigure[ ]{
    \begin{minipage}[t]{0.3\linewidth}
      \centering
      \includegraphics[width=1\linewidth]{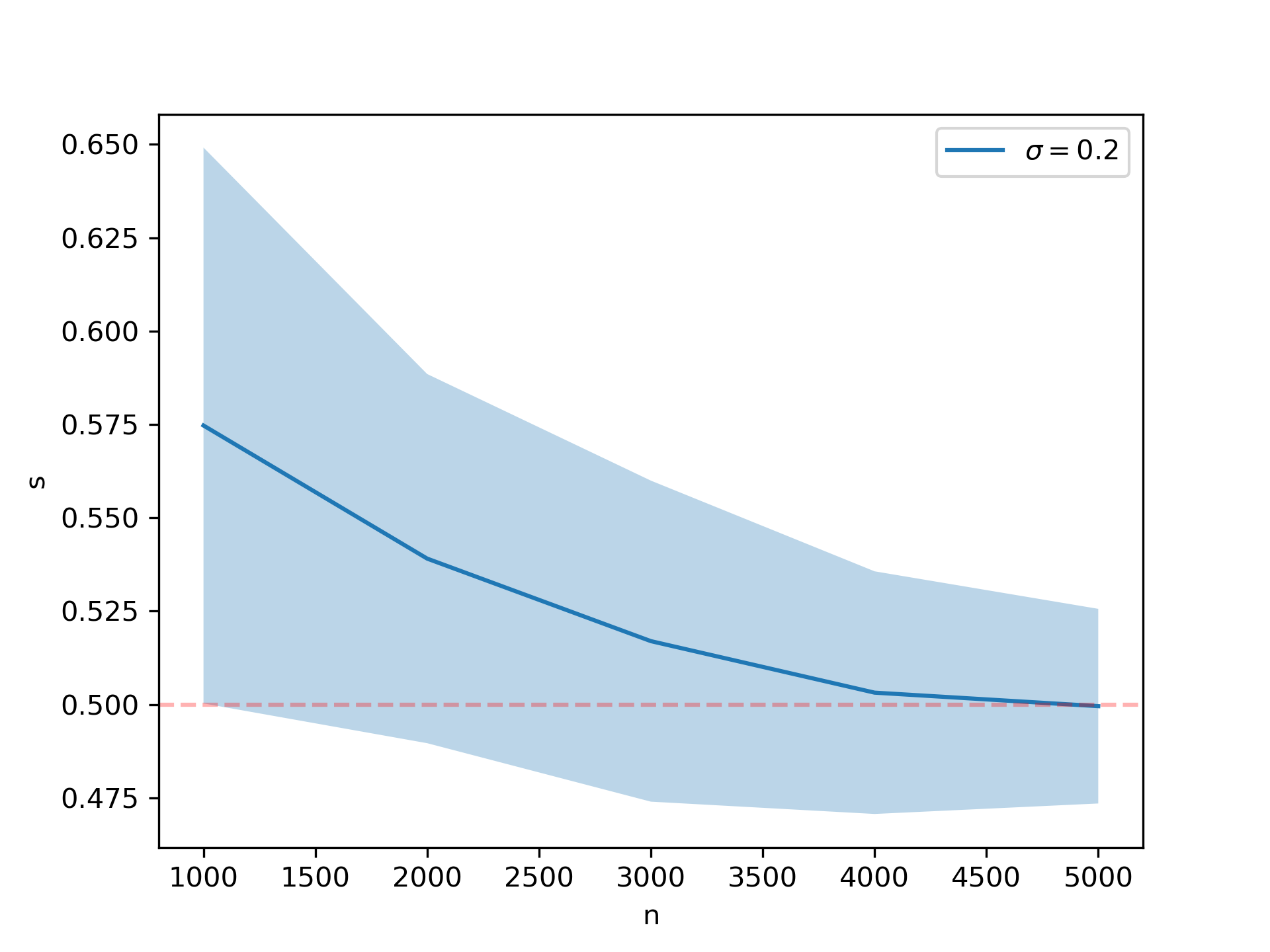}
    \end{minipage}
    }
  \centering
  \caption{Experiments for estimating the smoothness parameter $s$ in regression settings. (a) Naive estimation based on $2,000$ sample points for $\sigma=0$ (blue) and $\sigma=0.1$ (orange). (b) Truncation Estimation based on $2,000$ sample points with truncation point $100$. 
  In both plots (a) and (b), the $x$-axis is the logarithmic index $j$ and the $y$-axis is the logarithmic $p_j$.
  (c) Truncation Estimation across various values of sample size $n$, 
  each repeated 50 times. The blue line represents the average of estimates, the shaded area represents one standard deviation, and the true value is indicated by the orange dashed line. 
  }

  \label{fig:estimation_s_regression}
  \end{figure}
To better understand the estimation process, we first consider regression settings where the noises have an explicit form and we then consider classification settings. Suppose that we have $n$ i.i.d. samples of $X_{[n]}=[X_1,...,X_n]^{\top} $ and $Y_{[n]}=[Y_1,...,Y_n]^{\top} \in $ from $Y_i=f^*(X_i)+\sigma \epsilon_i$, where $\epsilon_i \sim \mathcal{N}(0,1)$. 

We start with a naive estimation method. Let $K_{min}(X_{[n]},X_{[n]})$ be the kernel matrix. Suppose the eigendecomposition is given by $K_{min}(X_{[n]},X_{[n]})= V\Sigma V^{\top}$, where $V = [v_1,...,v_n]$ is the eigenvector matrix, $v_i$'s are the eigenvectors, and $\Sigma$ is the diagonal matrix of the eigenvalues. We can estimate $r$ by estimating the decay rate of $p_j$, where $p_j = Y_{[n]}^{\top} v_j$. To visualize the convergence rate $r$, we perform logarithmic least-squares to fit $p_j$ with respect to the index $j$ and display the values of the slope $r$ and the smoothness parameter $s$.

For $\sigma=0$, $r$ can be accurately estimated by the above naive method since there is no noise in $Y_i$'s. The blue line and dots in Figure \ref{fig:estimation_s_regression} (a) present the estimation of $s$ in this case, where the estimate is around the true value $0.5$. However, for $\sigma=0.1$, the naive estimation is not accurate, as shown by the orange line and dots in Figure \ref{fig:estimation_s_regression} (a).

\begin{figure}[bhtp]
  \centering
\subfigure[ ]{
    \begin{minipage}[t]{0.33\linewidth}
      \centering
      \includegraphics[width=1\linewidth]{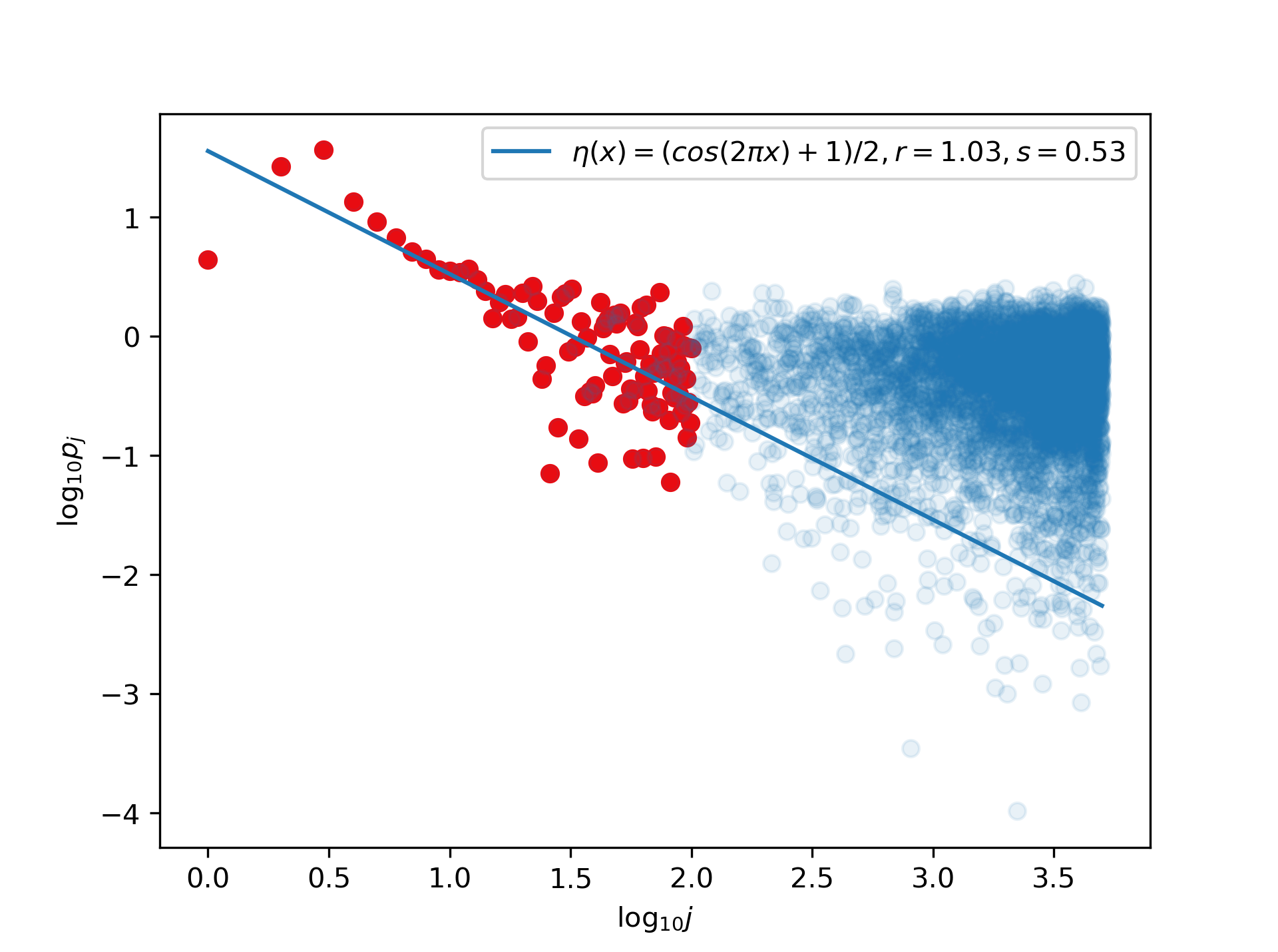}
    \end{minipage}%
    }
\subfigure[ ]{
    \begin{minipage}[t]{0.33\linewidth}
      \centering
      \includegraphics[width=1\linewidth]{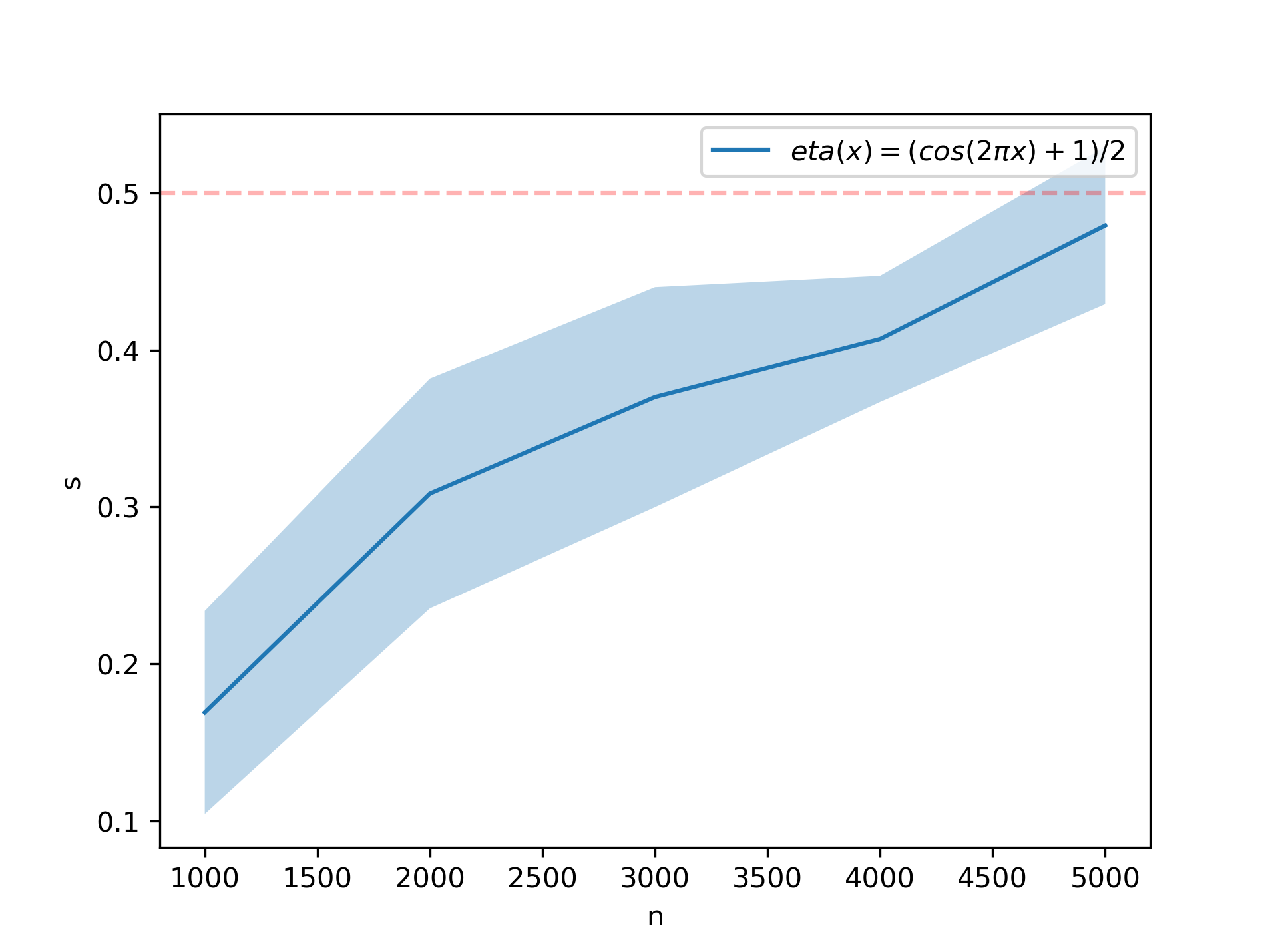}
    \end{minipage}%
    }
  \centering
  \caption{ Experiments for estimating the smoothness parameter $s$ in classification settings. (a) The experiment uses $5,000$ sample points and the truncation point is 100. (b)  Truncation Estimation across various values of sample size $n$, 
  each repeated 50 times. The blue line represents the average of estimates, the shaded area represents one standard deviation, and the true value is indicated by the orange dashed line. 
  }

  \label{fig:estimation_s_classification}
  \end{figure}

To improve the accuracy of the estimation, we introduce a simple modification called \textit{Truncation Estimation}, described as follows. We select some fixed integer as a truncation point and estimate the decay rate of $p_j$ up to the truncation point. 
For the example with $\sigma=0.1$, we choose the truncation point $100$ and the result is shown in Figure \ref{fig:estimation_s_regression} (b). We observe that the estimation becomes much more accurate than the naive estimation, with an estimate of $s=0.53$ not too far away from the true value $0.5$. 
In general, noise in the data can worsen the estimation accuracy, while increasing the sample size can improve the accuracy and robustness of the estimation. In Figure \ref{fig:estimation_s_regression} (c), we show the result for estimating $s$ in repeated experiments with more noisy data ($\sigma=0.2$), where we observe that as the sample size $n$ increases, the estimation becomes accurate. 
\paragraph{Estimation of $s$ in classification. }
Now we consider the classification settings, where the population is given by $\mathbb{P}(Y=1|X= x) = (f^*(x)+1)/2$. 
Unlike regression problems, the variance of the noise $\epsilon = y- f^*(x)$ is determined by $f^*(x)$ and may not be negligible. 
Nonetheless, in classification problems, we can still estimate the smoothness parameter $s$ using Truncation Estimation, thanks to the fact that increasing the sample size can improve its performance. 
The results are shown in Figure \ref{fig:estimation_s_classification}, where we can indeed make similar observations to those in Figure \ref{fig:estimation_s_regression} (b) and (c). 
  
\begin{table}[hbtp]
\centering
\begin{tabular}{p{1.5cm}| p{3cm}| p{3cm}| p{3cm}} 
 Kernel & MNIST& Fashion-MNIST & CIFAR-10\\
 \hline
NTK-1  & 0.4862 (0.0824) & 0.4417 (0.0934)& 0.1992 (0.0724)\\ [0.7ex] 
 \hline
NTK-2 &  0.4871 (0.0793)  & 0.4326 (0.0875)  & 0.2047 (0.0831) \\[0.7ex] 
 \hline
 NTK-3 &  0.4865 (0.0815) & 0.4372 (0.0768)  & 0.1965 (0.0795)  \\[0.7ex] 
 \hline

\end{tabular}
 \caption{Truncation Estimation of the relative smoothness $s$ of different real data sets with different NTKs. $NTK-L$ indicates the $L$-hidden-layer NTK. We only consider two classes of labels for each dataset: Label 1 and 7 for MNIST, trousers and sneakers for Fashion-MNIST, cars and horses for CIFAR-10. We randomly select 5,000 data points 
 and choose the truncation point $100$
 to estimate $s$. For each dataset and each kernel, we repeat 50 times and the standard deviation is in parentheses.}
\label{table}
\end{table}

As an application of Truncation Estimation, we estimate the relative smoothness of real data sets with respect to the NTK defined in \eqref{eq:NTK_Formula}. The results are shown in Table \ref{table}. We can see that with respect to the NTK, MNIST has the largest relative smoothness while CIFAR-10 has the smallest one. 
This result aligns with the common knowledge that MNIST is the easiest dataset while CIFAR-10 is the most difficult one of these three datasets. 

\paragraph{Limitations} The misspecified spectral algorithms (assuming $f^*_{\rho}\in[\mathcal{H}]^s$) are studied since 2009 (e.g., \cite{steinwart2009optimal, dicker2017kernel,PillaudVivien2018StatisticalOO,fischer2020sobolev,zhang2023optimality}). However, to the best of our knowledge, there is barely any result on the estimation of the smoothness $s$. This paper is the first to propose the $s$ estimation method even though the method is more susceptible to noise when the sample size is not enough or $f^*$ has more complex structures. For example, if $f^*= \sum_{j=1}^{\infty}f_je_j(x)$, where $f_j^2 = j^{-s_1\beta-1}$ when $j$ is odd and $f_j^2 = j^{-s_2\beta-1}$ when $j$ is even ($s_1>s_2$). For the kernel $K$ with EDR $\beta$, $f^*_\rho\in[\mathcal{H}]^{s_2}$ instead of $[\mathcal{H}]^{s_1}$ or $[\mathcal{H}]^{s}$ for some $s\in(s_2,s_1)$. In this mixed smoothness case, our method tends to give an estimation $\hat{s}\in(s_2,s_1)$. A more detailed discussion of the limitations is presented in the appendix. We will try to find more accurate $s$ estimation methods for general situations in the near future.

\section{Discussion}

In this paper, we study the generalization error of kernel classifiers in Sobolev space (the interpolation of the Sobolev RKHS). 
We show the optimality of kernel classifiers under the assumption that the ground true function is in the interpolation of RKHS with the kernel. 
The minimax optimal rate is $n^{-s\beta/2(s\beta+1)}$, where $s$ is the smoothness parameter of the ground true function.   
%
Building upon the connection between kernel methods and neural networks, we obtain an upper bound on the generalization error of overparameterized neural network classifiers. 
To make our theoretical result more applicable to real problems, we propose a simple method called Truncation Estimation to estimate the relative smoothness $s$. 
Using this method, we examine the relative smoothness of three real datasets, including MNIST, Fashion-MNIST and CIFAR-10. Our results confirm that among these three datasets, MNIST is the simplest for classification using NTK classifiers while CIFAR-10 is the hardest.  

\bibliography{iclr2024_conference}
\bibliographystyle{iclr2024_conference}

\section{Appendix}

In this section, we first show the proof of the upper bound of the classification excess risk (A.1 and A.2) and then present the minimax lower bound (A.3). Before the proof, We list again the standard assumptions for general RKHS $\mathcal{H}$ in this section.

\begin{assumption}[Source condition]\label{ass:rkhs}
  For $s > 0 $, there is a constant $B > 0 $ such that $f_{\rho}^{*} \in [\mathcal{H}]^{s}$ and
  \begin{displaymath}
    \| f_{\rho}^{*} \|_{[\mathcal{H}]^{s}} \le B.
  \end{displaymath}
\end{assumption}

\begin{assumption}[Eigenvalue Decay Rate (EDR)]\label{ass:eigendecay_rate}
    The EDR of the eigenvalues $\{\lambda_j\}$ associated to the kernel $K$ is $\beta$, i.e.,
    \begin{align}
        cj^{-\beta} \leq \lambda_j \leq Cj^{-\beta}
    \end{align}
    for some positive constants $c$ and $C$ and $\beta>1$. 
\end{assumption}

\begin{assumption}[Embedding index]\label{ass:emb}
Suppose that there exists $\alpha_{0} > 0$, such that
\begin{displaymath}
  \alpha_{0} = \inf\left\{ \alpha \in [\frac{1}{\beta},1] :  \left\|[\mathcal{H}]^\alpha \hookrightarrow L^{\infty}(\mathcal{X},\mu)\right\| < \infty  \right\},
\end{displaymath}
and we refer to $\alpha_{0}$ as the \textit{embedding index} of an RKHS $\mathcal{H}$.
\end{assumption}

Define the sampling operator $ K_{x}: \mathbb{R} \rightarrow \mathcal{H}, ~ y \mapsto y K(x,\cdot) $ and its adjoint operator $K_{x}^{*}: \mathcal{H} \rightarrow \mathbb{R},~ f \mapsto f(x)$. Further, we define the sample covariance operator  $T_X: \mathcal{H} \to \mathcal{H}$ as 
$$
T_X:=\frac{1}{n} \sum_{i=1}^n K_{X_i} K_{X_i}^*.
$$

Then we know that $\left\|T_X\right\| \leq\left\|T_X\right\|_1 \leq \kappa^2$, where $\|\cdot\|$ denotes the operator norm and $\|\cdot\|_1$ denotes the trace norm. Further, define the sample basis function 
$$
g_Z:=\frac{1}{n} \sum_{i=1}^n K_{X_i} Y_i \in \mathcal{H}.
$$

 We also introduce a more general framework known as spectra algorithm \cite{rosasco2005spectral,caponnetto2006optimal,bauer2007regularization}. We define the filter function and the spectral algorithms as follows:

\begin{definition}[Filter function] 
Let $\left\{\varphi_\nu:\left[0, \kappa^2\right] \rightarrow R ^{+} \mid \nu \in \Gamma \subseteq R ^{+}\right\}$ be a class of functions and $\psi_\nu(z)=1-z \varphi_\nu(z)$. If $\varphi_\nu$ and $\psi_\nu$ satisfy:
\begin{itemize}
    \item  $\forall \alpha \in[0,1]$, we have
\begin{equation}\label{filter1}
    \sup _{z \in\left[0, \kappa^2\right]} z^\alpha \varphi_\nu(z) \leq E \nu^{1-\alpha}, \quad \forall \nu \in \Gamma ;
\end{equation}

\item $\exists \tau \geq 1$ s.t. $\forall \alpha \in[0, \tau]$, we have
\begin{equation}\label{filter2}
    \sup _{z \in\left[0, \kappa^2\right]}\left|\psi_\nu(z)\right| z^\alpha \leq F_\tau \nu^{-\alpha}, \quad \forall \nu \in \Gamma,
\end{equation}
\end{itemize}

where $E, F_\tau$ are absolute constants, then we call $\varphi_\nu$ a filter function. We refer to $\nu$ as the regularization parameter and $\tau$ as the qualification.
\end{definition}

\begin{definition}[spectral algorithm] Let $\varphi_\nu$ be a filter function index with $\nu>0$. Given the samples $Z$, the spectral algorithm produces an estimator of $f_\rho^*$ given by
\begin{equation}\label{fv}
    \hat{f}_\nu=\varphi_\nu\left(T_X\right) g_Z.
\end{equation}
\end{definition}

\subsection{Some bounds}
Throughout the proof, we denote
$$
T_\nu=T+\nu^{-1} ; \quad T_{X \nu}=T_X+\nu^{-1}
$$
where $\nu$ is the regularization parameter. In addition, we denote $L^2( X , \mu)$ as $L^2, L^{\infty}( X , \mu)$ as $L^{\infty}$ for brevity throughout the proof. We use $a_n \asymp b_n$ to denote that there exist constants $c$ and $C$ such that $c a_n \leq b_n \leq C a_n, \forall n=1,2, \cdots ;$ use $a_n \lesssim b_n$ to denote that there exists an constant $C$ such that $a_n \leq C b_n, \forall n=1,2, \cdots$
In addition, denote the effective dimension as
$$
N (\nu)=\operatorname{tr}\left(T\left(T+\nu^{-1}\right)^{-1}\right)=\sum_{i \in N} \frac{\lambda_i}{\lambda_i+\nu^{-1}}
$$

\begin{lemma}\label{N nu}
Suppose $\nu>1$. 
     If $\lambda_i \asymp i^{-\beta}$, we have
$$
N (\nu) \asymp \nu^{\frac{1}{\beta}} .
$$
\end{lemma}
\begin{proof}
Since $c i^{-\beta} \leq \lambda_i \leq C i^{-\beta}$, we have
$$
\begin{aligned}
N (\nu) & =\sum_{i=1}^{\infty} \frac{\lambda_i}{\lambda_i+\nu^{-1}} \leq \sum_{i=1}^{\infty} \frac{C i^{-\beta}}{C i^{-\beta}+\nu^{-1}}=\sum_{i=1}^{\infty} \frac{C}{C+\nu^{-1} i^\beta} \\
& \leq \int_0^{\infty} \frac{C}{\nu^{-1} x^\beta+C} d x=\nu^{\frac{1}{\beta}} \int_0^{\infty} \frac{C}{y^\beta+C} d y \leq C_1 \nu^{\frac{1}{\beta}}
\end{aligned}
$$
for some constant $C_1$. Since $\nu>1$, the proof for the lower bound can be obtained similarly. 


\end{proof}

\subsubsection{Approximation error}
Recall that we have defined the sample basis function $g_Z$ and the spectral algorithm $\hat{f}_\nu$. We also need the following notations: define the expectation of $g_Z$ as
$$
g= E g_Z=\int_{ \mathcal{X}} K(x, \cdot) f_\rho^*(x) d \mu(x)=S_k^* f_\rho^* \in \mathcal{H} ,
$$
and
$$
f_\nu=\varphi_\nu(T) g=\varphi_\nu(T) S_k^* f_\rho^*
$$
The following conclusion based on \cite{zhang2023optimality} bounds the $L^2$-norm of $f_\nu-f_\rho^*$ for spectral algorithm:
\begin{lemma}\label{approxi}
    Suppose that Assumption \ref{ass:rkhs} holds for $0<s \leq 2 \tau$. Then for any $\nu>0$ , we have
$$
\left\|f_\nu-f_\rho^*\right\|_{L^2} \leq F_\tau R \nu^{-\frac{s}{2}} .
$$
\end{lemma}
\begin{proof}
    Because $f^*_{\rho} \in [\mathcal{H}]^s$, we assume $f^*_{\rho} = L^{\frac{s}{2}}g_0$ for some $g_0 \in L^2(P)$, so that $\|g_0\|_{L^2} \le B$ by Assumption \ref{ass:rkhs}. By the definition of $f_\nu$ , we have:
    $$
    \begin{aligned}
\|f_{\nu}-f^*_{\rho}\|_{L^2} &= \|\varphi_\nu(T) S_k^{*} f_\rho^*-f^*_{\rho}\|_{L^2}\\
&= \|(\varphi_\nu(L) L-Id)L^{\frac{s}{2}}g_0\|_{L^2}\\
& \le \|(\psi_{\nu}(L)L^{\frac{s}{2}}g_0\|_{L^2}\\
& \le F_r B \nu^{-\frac{s}{2}}.
\end{aligned}
    $$
Where the second equality holds by the definition of natural embedding inclusion operator $S_k$, and $Id$ denotes identity mapping. The first inequality holds because of the definition of $\psi_{\nu}$ and the second inequality holds for \eqref{filter2}.
\end{proof}

\subsubsection{Estimation error}
We rewrite the estimation error as follows
\begin{equation}\label{estimation}
    \begin{aligned}
\left\|\hat{f}_\nu-f_\nu\right\|_{L^2} & =\left\|T^{\frac{1}{2}} \left(\hat{f}_\nu-f_\nu\right)\right\|_{\mathcal{H}} \\
& =\left\|T^{\frac{1}{2}} T_\nu^{-\frac{1}{2}} \cdot T_\nu^{\frac{1}{2}} T_{X \nu}^{-\frac{1}{2}} \cdot T_{X \nu}^{\frac{1}{2}}\left(\hat{f}_\nu-f_\nu\right)\right\|_{\mathcal{H}} \\
& \leq\left\|T^{\frac{1}{2}} T_\nu^{-\frac{1}{2}}\right\|_{ \mathscr{B}  (\mathcal{H})} \cdot\left\|T_\nu^{\frac{1}{2}} T_{X \nu}^{-\frac{1}{2}}\right\|_{ \mathscr{B}  ( \mathcal{H} )} \cdot\left\|T_{X\nu}^{\frac{1}{2}}\left(\hat{f}_\nu-f_\nu\right)\right\|_{ \mathcal{H}} .
\end{aligned}
\end{equation}

\paragraph{Step 1.} 
The first part can be bounded by the following lemma, whose proof is simple and omitted. 
\begin{lemma}\label{1part}
	$$
    	\left\|T^{\frac{1}{2}}\left(T_\nu\right)^{-1 / 2}\right\|^2=\sup _{i \geq 1} \frac{\lambda_i}{\lambda_i+\nu^{-1}} \leq 1.
	$$
\end{lemma}

For the second part, we recall a result \cite[ Lemma 11]{fischer2020sobolev}.
\begin{lemma}[\cite{fischer2020sobolev}]\label{par2}
    Suppose that the RKHS~~$\mathcal{H}$ has the embedding index $\alpha_0$. Then for any $\alpha_0<\alpha \leq 1$ and all $\delta \in(0,1)$, with probability at least $1-\delta$, we have
$$
\left\|T_\nu^{-\frac{1}{2}}\left(T-T_X\right) T_\nu^{-\frac{1}{2}}\right\| \leq \frac{4 M_\alpha^2 \nu^\alpha}{3 n} B+\sqrt{\frac{2 M_\alpha^2 \nu^\alpha}{n}B} ,
$$
where
$$
B=\ln \frac{4 N (\nu)\left(\|T\|+\nu^{-1}\right)}{\delta\|T\|}. 
$$
\end{lemma}
\begin{lemma}\label{2part}
If the sample size $n \ge 8M_\alpha^2 \nu^\alpha B$, we have:
$$
\left\|T_\nu^{\frac{1}{2}} T_{X \nu}^{-\frac{1}{2}}\right\|_{ \mathscr{B}  ( \mathcal{H} )} \le 3
$$
holds with probability at least $1-\delta$.
\end{lemma}
\begin{proof}
We calculate $T_\nu^{\frac{1}{2}} T_{X \nu}^{-\frac{1}{2}}$ directly, and combined with Lemma \ref{par2}
$$
\begin{aligned}
	T_{X\nu}= T_{X}+\nu^{-1} &= T_{X} - T+ T +\nu^{-1}\\
	& = (T_\nu)^{\frac{1}{2}}\left[Id-(T_\nu)^{-\frac{1}{2}}(T-T_X)(T_\nu)^{-\frac{1}{2}}\right](T_\nu)^{\frac{1}{2}},
\end{aligned}
$$
By Lemma \ref{par2} when $n \ge 8M_\alpha^2 \nu^\alpha B$, we have:
$$
\left\|T_\nu^{-\frac{1}{2}}\left(T-T_X\right) T_\nu^{-\frac{1}{2}}\right\| \leq \frac{4}{3} \cdot\frac{1}{8}+\sqrt{ \frac{1}{4}} \leq \frac{2}{3} 
$$
holds with probability at least $1-\delta$. So the Neumann series gives us the following bound

 $$
 \begin{aligned}
\left\|T_\nu^{\frac{1}{2}} T_{X \nu}^{-\frac{1}{2}}\right\|_{ \mathscr{B}  ( \mathcal{H} )} &= \left\|\left(Id-(T_\nu)^{-\frac{1}{2}}(T-T_X)(T_\nu)^{-\frac{1}{2}}\right)^{-1}\right\|_{ \mathscr{B}  ( \mathcal{H} )}\\
&\le\sum_{k=0}^{\infty}\left\|T_\nu^{-\frac{1}{2}}\left(T-T_X\right) T_\nu^{-\frac{1}{2}}\right\|^k \\
&\le \sum_{k=0}^{\infty}\left(\frac{2}{3}\right)^k= 3.
\end{aligned}
 $$
\end{proof}
\paragraph{Step 2.} For the third part in the last line of \eqref{estimation}, we have
\begin{equation}\label{3part}
    \left\|T_{X \nu}^{\frac{1}{2}}\left(\hat{f}_\nu-f_\nu\right)\right\|_{ \mathcal{H} } \leq\left\|T_{X \nu}^{\frac{1}{2}} \varphi_\nu\left(T_X\right)\left(g_Z-T_X f_\nu\right)\right\|_{ \mathcal{H} }+\left\|T_{X \nu}^{\frac{1}{2}} \psi_\nu\left(T_X\right) f_\nu\right\|_{ \mathcal{H} } .
\end{equation}

The second part of RHS in \eqref{3part} is complicated to calculate, but its proof follows  the same argument as in Step 3 of Theorem 16 in \cite{zhang2023optimality}. 
Therefore, we provide the following result without proof: 
\begin{lemma}[Theorem 16 in \cite{zhang2023optimality}]
    If $s<2\tau$, we have:
\begin{equation}\label{3part2}
    \left\|T_{X \nu}^{\frac{1}{2}} \psi_\nu\left(T_X\right) f_\nu\right\|_{ \mathcal{H} } \leq 6 F_\tau E R \nu^{-\frac{s}{2}}+\Delta_1 I_{s>2},
\end{equation}

where $\Delta_1$ denotes
$$
\Delta_1  = 32 \max \left(\frac{s-1}{2}, 1\right) E F_\tau R \kappa^{s-1} \nu^{-\frac{1}{2}} n^{-\frac{\min (s, 3)-1}{4}} \ln \frac{6}{\delta}.
$$
We need $s < 2\tau$, because we use \eqref{filter2} to upper bound $\|T_X^{\frac{s}{2}}\psi_{\nu}(T_X)\|$ for every $s$.
\end{lemma}

To bound the first term of RHS in \eqref{3part}, we begin with a lemma, whose proof is postponed to Section \ref{pf bernstein}. 
\begin{lemma}\label{bernstein}
Suppose that Assumption \ref{ass:rkhs}, \ref{ass:eigendecay_rate} and \ref{ass:emb} hold for $\frac{1}{\beta}\le \alpha_0 < 1$, and conditional probability of $Y$ is given by 
$$
\mathbb{P}(Y =1|X = x) = \frac{1+f_\rho^*(x)}{2}.
$$
Then given $\nu>0$, and $n \geq 1$, for any $\delta>0$ and $\alpha> \alpha_0$,  the following bound is satisfied with probability not less than $1-\delta$
	$$
	\begin{aligned}
		 \left\|\left(T_\nu\right)^{-1 / 2}\left(\left(g_Z-T_X f_{\nu}\right)-\left(g -T f_{\nu}\right)\right)\right\|_{\mathcal{H}}^2 
		\leq  \frac{128 \left(\log\frac{2}{\delta}\right)^2}{n}\left( N(\nu)+ \frac{4M_\alpha^2\nu^\alpha}{n }\right).
	\end{aligned}
	$$
\end{lemma}

The next lemma provides a bound on the first term of RHS in \eqref{3part}. 

\begin{lemma}\label{3part3}
Suppose that Assumption \ref{ass:rkhs}, \ref{ass:eigendecay_rate} and \ref{ass:emb} hold for $\frac{1}{\beta}\le \alpha_0 < 1$. If the sample size $n \geq 8M_\alpha^2 \nu^\alpha B$, where B is defined in Lemma \ref{par2} and  $\nu>0$, then for any $\delta>0$ and $\alpha> \alpha_0$,  the following bound is satisfied with probability at least $1-\delta$:
    $$
    \left\|T_{X \nu}^{\frac{1}{2}} \varphi_\nu\left(T_X\right)\left(g_Z-T_X f_\nu\right)\right\|_{ \mathcal{H}} \le C\left[\left(\log\frac{2}{\delta}\right)\left( \frac{N(\nu)^{\frac{1}{2}}}{\sqrt{n}}+ \frac{2M_\alpha\nu^\frac{\alpha}{2}}{n }\right)+\left\|f_\rho^*-f_\nu\right\|_{L^2} \right],
    $$
    where C is an absolute constant.
\end{lemma}
\begin{proof}
    \begin{equation}
    \begin{aligned}
\left\|T_{X \nu}^{\frac{1}{2}} \varphi_\nu\left(T_X\right)\left(g_Z-T_X f_\nu\right)\right\|_{ \mathcal{H}} & =\left\|T_{X \nu}^{\frac{1}{2}} \varphi_\nu\left(T_X\right) T_{X \nu}^{\frac{1}{2}} \cdot T_{X \nu}^{-\frac{1}{2}} T_\nu^{\frac{1}{2}} \cdot T_\nu^{-\frac{1}{2}}\left(g_Z-T_X f_\nu\right)\right\|_{  \mathcal{H} } \\
& \leq\left\|T_{X \nu}^{\frac{1}{2}} \varphi_\nu\left(T_X\right) T_{X \nu}^{\frac{1}{2}}\right\|_{ \mathscr{B} ( \mathcal{H} )} \cdot\left\|T_{X \nu}^{-\frac{1}{2}} T_\nu^{\frac{1}{2}}\right\|_{ \mathscr{B} (  \mathcal{H} )}\\
&\cdot\left\|T_\nu^{-\frac{1}{2}}\left(g_Z-T_X f_\nu\right)\right\|_{  \mathcal{H}} .
\end{aligned}
\end{equation}
 
The first term can be upper bounded as
$$
\begin{aligned}
\left\|T_{X \nu}^{\frac{1}{2}} \varphi_\nu\left(T_X\right) T_{X \nu}^{\frac{1}{2}}\right\|_{ \mathscr{B} ( \mathcal{H} )}& = \|(T_X+\nu^{-1})\varphi_\nu\left(T_X\right)\|_{ \mathscr{B} ( \mathcal{H} )} \\
&\le \|T_X\varphi_\nu\left(T_X\right)\|_{ \mathscr{B} ( \mathcal{H})}+\nu^{-1}\|\varphi_\nu\left(T_X\right)\|_{ \mathscr{B} ( \mathcal{H})}\\
&\le 2E 
\end{aligned}
$$
because we have $ z \varphi_\nu(z) \leq E $ by \eqref{filter1} and $\varphi_\nu(z) \leq E \nu $ by \eqref{filter2}. 
Using Lemma \ref{2part}, the second term can be bounded as
 $$
\left\|T_{X \nu}^{-\frac{1}{2}} T_\nu^{\frac{1}{2}}\right\|_{ \mathscr{B} ( \mathcal{H} )}  \le 3.
$$
It remains to bound the third term. Lemma \ref{bernstein} shows that 
\begin{align}
		& \left\|\left(T_\nu\right)^{-1 / 2}\left(\left(g_Z-T_X f_{\nu}\right)-\left(g -T f_{\nu}\right)\right)\right\|_{\mathcal{H}}^2 
		\leq  \frac{128 \left(\log\frac{2}{\delta}\right)^2}{n}\left( N(\nu)+ \frac{4M_\alpha^2\nu^\alpha}{n }\right).
\end{align}
Thus, we have
$$
\begin{aligned}
\left\|T_\nu^{-\frac{1}{2}}\left(g_Z-T_X f_\nu\right)\right\|_{ \mathcal{H} } & \leq\left\|T_\nu^{-\frac{1}{2}}\left[\left(g_Z-T_X f_\nu\right)-\left(g-T f_\nu\right)\right]\right\|_{ \mathcal{H} }+\left\|T_\nu^{-\frac{1}{2}}\left(g-T f_\nu\right)\right\|_{ \mathcal{H} } \\
& \leq {8\sqrt{2} \left(\log\frac{2}{\delta}\right)}\left( \frac{N(\nu)^{\frac{1}{2}}}{\sqrt{n}}+ \frac{2M_\alpha\nu^\frac{\alpha}{2}}{n }\right)+\left\|T_\nu^{-\frac{1}{2}} S_k^*\right\|_{ \mathscr{B} \left(L^2, \mathcal{H} \right)}\left\|f_\rho^*-f_\nu\right\|_{L^2} \\
& \leq {8\sqrt{2} \left(\log\frac{2}{\delta}\right)}\left( \frac{N(\nu)^{\frac{1}{2}}}{\sqrt{n}}+ \frac{2M_\alpha\nu^\frac{\alpha}{2}}{n }\right)+\left\|f_\rho^*-f_\nu\right\|_{L^2},
\end{aligned}
$$
where the second term is approximation error that has been referred to in Lemma \ref{approxi}.
\end{proof}

\paragraph{Step 3.} 
 Now we combine the bounds for the three parts of estimation error  in \eqref{estimation}: the first two parts are corresponding to Lemma \ref{1part} and \ref{2part} respectively, and the third part is corresponding to \eqref{3part}, \eqref{3part2}, and Lemma \ref{3part3}. Then  based on Assumptions the same as  Lemma \ref{3part3}, combined with $s< 2\tau$, we conclude that for any $\delta>0$, it holds with probability at least $1-\delta$ 
\begin{equation}\label{estimation2}
\begin{aligned}
    \left\|\hat{f}_\nu-f_\nu\right\|_{L^2} &\le C\left(\log\frac{2}{\delta}\right)\left( \frac{N(\nu)^{\frac{1}{2}}}{\sqrt{n}}+ \frac{2M_\alpha\nu^\frac{\alpha}{2}}{n }\right)\\
    &+\left\|f_\rho^*-f_\nu\right\|_{L^2}+F_\tau E R \nu^{-\frac{s}{2}}+\Delta_1 I_{s>2}, 
    \end{aligned}
\end{equation}
where C is an absolute constant.

\subsubsection{Proof of Lemma \ref{bernstein}}\label{pf bernstein}

\begin{lemma}[Lemma 13 in \cite{fischer2020sobolev}]\label{sigma L}
	Let $(\mathcal{X}, B )$ be a measurable space, $\mathcal{H}$ be a separable RKHS on $\mathcal{X}$ w.r.t. a bounded and measurable kernel $k$, and $\mu$ be a probability distribution on $\mathcal{X}$. Then the following equality is satisfied, for $\nu>0$,
	$$
	\int_{\mathcal{X}}\left\|\left(T_\nu\right)^{-1 / 2} K(x, \cdot)\right\|_{\mathcal{H}}^2 d \mu(x)= N (\nu) .
	$$
	If, in addition, $M_\alpha<\infty$ is satisfied, then the following inequality is satisfied, for $\nu>0$ and $\mu$-almost all $x \in \mathcal{X}$,
	$$
	\left\|\left(T_\nu\right)^{-1 / 2} K(x, \cdot)\right\|_{\mathcal{H}}^2 \leq M_\alpha^2 \nu^{-\alpha}.
	$$

we also consider $C_x: \mathcal{H} \rightarrow \mathcal{H}$ the integral operator w.r.t. the point measure at $x \in \mathcal{X}$,
$$
C_x f:=f(x) K(x, \cdot)=\langle f, K(x, \cdot)\rangle_\mathcal{H} K(x, \cdot),
$$
And we have the operator norm:
	$$
	\left\|\left(T_\nu\right)^{-1 / 2}C_x\left(T_\nu\right)^{-1 / 2}\right\|_{\mathscr{B}(\mathcal{H})}=	\left\|\left(T_\nu\right)^{-1 / 2} K(x, \cdot)\right\|_\mathcal{H}^2\leq M_\alpha^2 \nu^{\alpha}.
	$$
\end{lemma}

\begin{lemma}[Bernstein's Inequality in \cite{caponnetto2007optimal}]\label{bern1}
	
 Let $(\Omega, B , P)$ be a probability space, $H$ be a separable Hilbert space, and $\xi: \Omega \rightarrow H$ be a random variable with
	$$
	\begin{aligned}
	E _P[\|\xi\|_{\mathcal{H}}^2] &\leq \sigma^2,\\
	\|\xi(\omega)\|_{\mathcal{H}} &\leq \frac{L}{2}, \quad a.s.
	\end{aligned}
	$$

 Then, for $\tau \geq 1$ and $n \geq 1$, the following concentration inequality is satisfied
	$$
	\mathbb{P}^n\left(\left(\omega_1, \ldots, \omega_n\right) \in \Omega^n:\left\|\frac{1}{n} \sum_{i=1}^n \xi\left(\omega_i\right)- E _P \xi\right\|_{\mathcal{H}}^2 \geq 32 \frac{\tau^2}{n}\left(\sigma^2+\frac{L^2}{n}\right)\right) \leq 2 e^{-\tau}.
	$$
	\end{lemma}
\begin{proof}[Proof of Lemma \ref{bernstein}]
	Consider the random variable $\xi: X \times \{-1,1\} \to \mathcal{H}$:
	$$
	\xi(x,y) = \left(T_\nu\right)^{-1 / 2}(y - f_{\nu}(x)) K(x, \cdot).
	$$
Denoted by $D$ the empirical measure corresponding to the sample $\{X_i,Y_i\}_{i=1}^n$. It holds that 
	$$
	\begin{aligned}
		\frac{1}{n} \sum_{i=1}^n\left(\xi\left(X_i, U_i\right)- E _P \xi\right) & = E _D \xi- E \xi \\
		& =\left(T_\nu\right)^{-1 / 2}\left(\left(g_Z-T_X f_{\nu}\right)-\left(g -T f_{\nu}\right)\right).
	\end{aligned}
	$$
	By Lemma \ref{sigma L}, it holds that 
	$$
	\|\xi(x,y)\|_{\mathcal{H}} \le 2	\left\|\left(T_\nu\right)^{-1 / 2} K(x, \cdot)\right\|_{\mathcal{H}} \le 2M_\alpha \nu^{\frac{\alpha}{2}}
	$$
for $\mu$-almost all $x\in \mathcal{X}$, and it holds that 
	$$
	\begin{aligned}
	E(\|\xi(x,y)\|_{\mathcal{H}}^2)  &= \int_{\mathcal{X}}\left\|\left(T_\nu\right)^{-1 / 2} K(x, \cdot)\right\|_{\mathcal{H}}^2\left[\sum_{i=1,-1}\left|i-f_{\nu}(x)\right|^2 P(Y=i\mid X=x)\right] d \mu(x)\\
	&\le 4N(\nu).
	\end{aligned}
	$$
 We let $\sigma^2 = 4N(\nu)$ and $L =  4M_\alpha \nu^{\frac{\alpha}{2}}$ and apply the Bernstein inequality in Lemma \ref{bern1} to complete the proof. 
\end{proof}

\subsection{Upper bound on excess risk}
\begin{theorem}[$L^2$-risk upper bound]\label{regression}
    Suppose that  Assumption \ref{ass:rkhs}, \ref{ass:eigendecay_rate} and \ref{ass:emb} holds, and that $0<s \leq 2 \tau$. Let $\hat{f}_\nu$ be the estimator defined by \eqref{fv}. Then by choosing $\nu \asymp n^{\frac{\beta}{s \beta+1}}$, for any fixed $\delta \in(0,1)$ and any $1\ge\alpha >\alpha_0$, when $n$ is sufficiently large, with probability at least $1-\delta$, we have
$$
\left\|\hat{f}_\nu-f_\rho^*\right\|_{L^2}^2 \leq\left(\ln \frac{6}{\delta}\right)^2 C n^{-\frac{s\beta}{s \beta+1}},
$$
where $C$ is a constant independent of $n$ and $\delta$.
\end{theorem}
\begin{proof}
We decomposed  $L^2$-risk into the sum of the approximation error and the estimation error as 
as $$\|\hat{f}_\nu-f_\rho^*\|_{L^2} \le  \|f_\nu-f_\rho^*\|_{L^2}+\|\hat{f}_\nu-f_\nu\|_{L^2}. $$
Using Lemma \ref{approxi} for the appoximation error and  \eqref{estimation2} for the estimation error, we have
$$
\begin{aligned}
    \|\hat{f}_\nu-f_\rho^*\|_{L^2} &\le C\left(\log\frac{2}{\delta}\right)\left( \frac{N(\nu)^{\frac{1}{2}}}{\sqrt{n}}+ \frac{2M_\alpha\nu^\frac{\alpha}{2}}{n }\right)+F_\tau (E+1) R \nu^{-\frac{s}{2}}+\Delta_1 I_{s>2}.  
    \end{aligned}
$$
Choosing $\nu \asymp n^{\frac{\beta}{s \beta+1}}$,  we can obtain the following rates:
\begin{itemize}
    \item By Lemma~\ref{N nu}, 
    \begin{equation}\label{th31}
\frac{N (\nu)^{\frac{1}{2}}}{\sqrt{n}} \asymp \frac{\nu^{\frac{1}{2 \beta}}}{\sqrt{n}}=n^{-\frac{1}{2} \frac{s \beta}{s \beta+1}}
\end{equation}
\item 
\begin{equation}\label{th32}
    F_\tau (E+1) R \nu^{-\frac{s}{2}} \asymp n^{-\frac{1}{2} \frac{s \beta}{s \beta+1}}
\end{equation}
\item Following Theorem 16 in \cite{zhang2023optimality}, for any $s> 2$, we  have
\begin{equation}\label{th33}
\Delta_1 \lesssim \left(\ln \frac{6}{\delta}\right)n^{-\frac{1}{2} \frac{s \beta}{s \beta+1}}
\end{equation}
\item Since $\alpha_0 = \frac{1}{\beta}$ for Sobolev RKHS, 
we can choose some $\alpha$ in $(\alpha_0 , s+\frac{2}{\beta})$. It follows that 
    \begin{equation}\label{th34}
\frac{\nu^{\frac{\alpha}{2}}}{n} \lesssim n^{-\frac{1}{2} \frac{s \beta}{s \beta+1}}
\end{equation}
\end{itemize}

Combining the above \eqref{th31}, \eqref{th32},\eqref{th33} and \eqref{th34}, we complete the proof.
\end{proof}

\begin{proof}[Proof of Theorem 
2
]
The classification excess risk can be rewritten as
    \begin{equation}\label{excessrisk}
    \begin{aligned}
    &\mathcal{E}(\hat{f}) = \mathcal{L}(\hat{f}) - \mathcal{L}^* \\
    = &\bm{E}_{X} [(1-\eta(X))(\mathbb{I}\{\hat{f}(X)\geq 0\} -  \mathbb{I}\{2\eta(X)-1\geq 0\})+ \eta(X) (\mathbb{I}\{\hat{f}(X)<0\} -\mathbb{I}\{2\eta(X)-1< 0\}))]\\
     =&\bm{E}_{X}(|2\eta(X)-1| \mathbb{I}\{\hat{f}(X)(2\eta(X)-1)<0\} )\\
    =&\bm{E}_{X}(|f^*_{\rho}(X)| \mathbb{I}\{\hat{f}(X)f^*_{\rho}(X)<0\} ).
\end{aligned}
\end{equation}
    
    Based on \eqref{excessrisk} and the classic upper classification upper bound by \cite{devroye2013probabilistic}, we have:
    $$
    \begin{aligned}
    \mathcal{E}(\hat{f}_{\nu}) &=\bm{E}_{X}(|f^*_\rho(X)| \mathbb{I}\{\hat{f}_\nu(X)f^*_\rho(X)<0\} )\\
    &\le \int_{\mathcal{X}}\left|f^*_\rho(X)-\hat{f}_\nu(X)\right| d\nu(X)\\
    & \le  \|\hat{f}_\nu-f_\rho^*\|_{L^2}
    \end{aligned}
    $$
   The remaining proof directly follows from  Theorem \ref{regression}.
\end{proof}

\subsection{Minimax lower bound}
\begin{proposition}[Theorem 2.5 in \cite{tsybakov2009_IntroductionNonparametric}]\label{prop:tsy}
    Assume that $M\ge 2$ and suppose that $\Theta$ contains elements $\theta_0,\theta_1,\dots, \theta_M$ and $P_{\theta_0},P_{\theta_1},\dots, P_{\theta_M}$ are the probability measures such that 
    \begin{itemize}
        \item[(i)] $d(\theta_i,\theta_j)\ge 2s>0$, $\forall 0\le i\le j \le M$;
        \item[(ii)] $P_{\theta_j}\ll P_{\theta_0}$, $\forall j=1,\dots,M$ and 
        \begin{align}
            \frac{1}{M}\sum_{j=1}^{M}KL(P_{\theta_j},P_{\theta_0})\le a \log M
        \end{align}
        with $0<\alpha<1/8$. 
    \end{itemize}
    Then
    \begin{align}
        \inf_{\hat{\theta}}\sup_{\theta\in\Theta} \mathbb{P}_{\theta}(d(\hat{\theta},\theta)>s)\ge \frac{\sqrt{M}}{1+\sqrt{M}}(1-2\alpha -\sqrt{\frac{2\alpha}{\log M}})>0
    \end{align}
\end{proposition}

\begin{lemma}[Varshamov-Gilbert Bound]\label{lem:VG_bound}
    Given $m\ge 8$, there exist $M\ge 2^{m/8}$ different elements on $\omega^{(0)},...,\omega^{(M)}$ on $\{-1,1\}^{m}$ and $\omega^{(0)}=(0,...,0)$ such that 
    \begin{align}
        \sum_{k=1}^{m}|\omega^{(i)}_k-\omega^{(j)}|\ge \frac{m}{4},\quad 0\le i<j\le M.
    \end{align}
\end{lemma}

\begin{lemma}
    For $r>\frac{d}{2}$, $H^{r}(\mathcal{X})$ is a separable RKHS with respect to a bound kernel and the corresponding EDR is 
    \begin{align*}
        \beta=\frac{2r}{d}.
    \end{align*}
\end{lemma}

    Let $u:\mathbb{R}_{+} \to \mathbb{R}_{+}$ be a nonincreasing infinitely differentiable function such that $u=1$ on $[0,1/4]$ and $u=0$ on $[1/2,\infty)$. We can take $u(x)=\left( \int_{1/4}^{1/2} u_1(z) \mathrm{d}z \right)^{-1} \int_{x}^{\infty} u_1(z) \mathrm{d}z$ where 
    \begin{align}
        u_1(x) = \begin{cases}
        \exp\left\{ -\frac{1}{(1/2-x)(x-1/4)}\right\}, & \text{for $x\in(1/4,1/2)$,}\\
        0,  & \text{otherwise}.
        \end{cases}
    \end{align}
    Given an integral $q=q(n)\ge 1$, we define the regular grid on $\mathbb{R}^d$ as
    \begin{align}
        G_q=\left\{\left(\frac{2k_1+1}{2q},...,\frac{2k_d+1}{2q}\right) : k_i\in\{0,...,q-1 \}, i=1,...,d \right\}.
    \end{align}
    We consider the partition $\mathcal{X}_1,...,\mathcal{X}_{q^d}$ of $[0,1]^d$ canonically defined using the grid $G_q$. $x\in\mathcal{X}_k$ if $z_k\in G_q$ is the closest point to $x\in [0,1]^{d}$. If there exist several points in $G_q$ closest to $x$ we define $x\in\mathcal{X}_k$ if $z_k$ is closest to $0$.

    \begin{lemma}\label{lem:phi_x_in_H_r}
     $\psi(x) = C_\psi q^{-sr}\sum_{k=1}^{q^{d}} \phi (q[x-z_k]) \in H^{sr}(\mathcal{X})$, where $\phi(x)=u(\|x\|_2)$.

    \end{lemma}

    \begin{proof}[Proof of Lemma \ref{lem:phi_x_in_H_r}]

    By the definition of $\phi$, we have $\phi\in H^{sr}$ for any fixed $sr>0$. Thus, $\|\phi\|_{H^{sr}}^2$ is bounded.
    
    Denote $\psi_k(x)=\phi (q[x-z_k])$ and $\varphi(x) = \psi_k(x+z_k) = \phi(qx)$. It is easy to find that: $i)$
    $\|\psi_k\|_{[H^{sr}]}=\|\varphi\|_{H^{sr}}$ ; $ii)$ $\langle \psi_i, \psi_j \rangle_{H^{sr}} =0$ for $i\neq j$. 
    \begin{align}
        \|\psi\|_{H^{sr}}^2 &= C_\psi^2 q^{-2sr} \|\sum_{k=1}^{q^{d}}\psi_k\|_{H^{sr}}^2\\
        &= C_\psi^2 q^{-2sr} \left(\sum_{k=1}^{q^{d}}\|\psi_k\|_{H^{sr}}^2 + 2\sum_{i\neq j} \langle \psi_i, \psi_j \rangle_{H^{sr}}^2 \right)\\
        &=C_\psi^2 q^{d-2sr}  \|\varphi\|_{H^{sr}}^2
    \end{align}

    Denote the Fourier transform of $\phi$ and $\psi$ as $\hat{\phi}$ and $\hat{\psi}$
    \begin{align}
        \hat{\varphi}(\xi) &=\int_{B(0,\frac{1}{2q})} \phi(qx) e^{-2\pi i \xi x}\mathrm{d}x\\
        &= q^{-d}\int_{B(0,\frac{1}{2})} \phi(y) e^{-2\pi i \frac{\xi}{q} y} \mathrm{d}y\\
        &=  q^{-d} \hat{\phi}(\frac{\xi}{q})
    \end{align}
    
    Since $u(x)$ is infinitely differentiable function on $[0,1]^d$, then $\|\phi \|_{H^{s}}$ is bounded for any fixed $s>0$. Then
        \begin{align}
            &\|\varphi\|_{H^{sr}}^2 =q^{-2d} \int_{R^{d}} |\hat{\phi}(\frac{\xi}{q})|^2 (1+\|\xi\|^2_2)^{s}\mathrm{d}\xi\\
            &=q^{-2d} \int_{R^{d}} |\hat{\phi}(z)|^2 (1+q^2\|z\|^2_2)^{sr} q^{d}\mathrm{d}z\\
            &\le q^{2sr-d} \int_{R^{d}} |\hat{\phi}(z)|^2 (1+\|z\|^2_2)^{sr} \mathrm{d}z\\
            &=q^{2sr-d} \|\phi\|_{H^{sr}}^2.
        \end{align}
    Thus, $\|\psi\|_{H^sr}^2 \leq C_\psi^2  \|\phi\|_{H^{sr}}^2$.
    \end{proof}

\begin{proof}[Proof of Theorem 
1
]
Denote $m= q^d$, we define $\mathcal{X}_0=\mathbb{R}^d \setminus \bigcup_{i=1}^{m}\mathcal{X}_i $. Thus, $\mathcal{X}_0,...,\mathcal{X}_m$ form a partition of $\mathbb{R}^{d}$.

Define the hypercube $\mathcal{C}=\{P_{\omega}:\omega=(\omega_1,...,\omega_m )\in \{-1,1\}^{m}\}$ of probability distribution $P_{\omega}$ of $(X,Y)$ on $Z=\mathbb{R}^{d}\times \{0,1\}$ as follows.

For any $P_{\omega}$, the marginal distribution of $X$ does not depend on $\omega$ and has a density $\mu$ w.r.t. the Lebesgue measure on $\mathbb{R}^{d}$ defined in the following way. Denote $v= m^{-1}$. Let $\mu(x)=v/\lambda[B(0,(4q)^{-1})]$ for $x\in B(z,(4q)^{-1})$, $z\in G_q$ and $\mu(x)=0$ otherwise.

By Lemma \ref{lem:VG_bound}, there exist $M\ge 2^{m/8}$ different elements on $\omega^{(0)},...,\omega^{(M)}$ on $\{-1,1\}^{m}$ and $\omega^{(0)}=(0,...,0)$. We take 
\begin{align}
    f_i(x)=\omega^{(i)}_k\psi(x), \quad x\in \mathcal{X}_k, \quad i\in[M]
\end{align}
and $f_i(x)=0$ for $x\in \mathcal{X}_0$. $f_0=0$ for $x\in\mathcal{X}$. We will assume that $C_\psi\leq1$ to ensure that $\eta_{i}(x)=(1+f_i(x))/2$ take values in $[0, 1]$.

By Lemma \ref{lem:phi_x_in_H_r}, we have $f_i\in H^{sr}= [\mathcal{H}]^{s}$, $i=0,...,M$.
Since $\eta_{i}(x)=(1+f_i(x))/2$, we have 
\begin{align}
    &KL(\rho_i^n,\rho_0^n) = n KL(\rho_i,\rho_0) \\
    &= n \sum_{k=1}^{m} \int_{\mathcal{X}_k} \left(\frac{1+\omega_{k}^{(i)}\psi(x)}{2} \ln( 1+\omega_{k}^{(i)}\psi(x) ) + \frac{1-\omega_{k}^{(i)}\psi(x)}{2} \ln( 1-\omega_{k}^{(i)}\psi(x)) \right)\mu(x) \mathrm{d}x\\
    &\le \frac{1}{2} n m  \int_{B(x_0,(4q)^{-1})} \frac{v}{\lambda[B(0,(4q)^{-1})]}\left( \ln(1-\psi^2(x)) + \psi(x) \ln\left( \frac{1+\psi(x)}{1-\psi(x)} \right) \right) \mathrm{d}x\\
    &\le \frac{1}{2} n m  \int_{B(x_0,(4q)^{-1})} \frac{v}{\lambda[B(0,(4q)^{-1})]}\left(  \psi(x) \ln\left( \frac{1+\psi(x)}{1-\psi(x)} \right) \right) \mathrm{d}x\\
    &\le C n m v q^{-2sr}
\end{align}
The last inequality is because with sufficiently large $n$,  $\psi(x)=C_\psi q^{-sr}<1/2$ for $x \in B(x_0,(4q)^{-1})$. To satisfy the second condition of Proposition \ref{prop:tsy}, we need $C_\psi n m v q^{-2sr} \leq a \frac{m}{8}\ln2  $. We can take $n v q^{-2sr} = \Theta(1)$. Thus, $q=n^{1/(2sr+d)}$.

For the first condition, we have
\begin{align}
    &d(f_i,f_j) = \mathbb{E}_x [|f_i(x)|1_{f_i(x)f_j(x)<0}]\\
    &= \mathbb{E}_x \left[\sum_{k=1}^{m} \psi(x) 1_{\sigma_k^{(i)}\sigma_k^{(j)}<0} 1_{x\in \mathcal{X}_k} \right]\\
    &\geq \frac{C_\psi q^{-sr} m v}{4}\\
    &= C n^{-sr/(2sr+d)}
\end{align}

for some constant $C>0$. By Proposition \ref{prop:tsy}, we have the minimax rate $n^{-sr/(2sr+d)}$. Since $r=\beta d/2$, we have the minimax rate $n^{-\frac{s\beta}{2(s\beta+1)}}$.
\end{proof}

\subsection{Embedding index of Sobolev and dot-product kernels}

\subsubsection{Sobolev kernel}

The interpolation space of $H^{r}(\mathcal{X})$ under Lebesgue measure is given by
\begin{equation}\label{sobolev_interpolation}
  [H^{r}(\mathcal{X})]^{s} = H^{rs}(\mathcal{X}).
\end{equation}

By the embedding theorem of (fractional) Sobolev space (Theorem 4.27 in \cite{adams2003sobolev}), letting $\theta = r-\frac{d}{2}>0$, we have
$$
H^r( \mathcal{X} ) \hookrightarrow C^{0, \theta}( \mathcal{X} )\hookrightarrow L^{\infty}( \mathcal{X} ), \quad \theta=r-\frac{d}{2}.
$$
Combined with \eqref{sobolev_interpolation}, for a Sobolev RKHS $\mathcal{H} =H^r( \mathcal{X} ), r>\frac{d}{2}$ and any $\alpha>\frac{1}{\beta}=\frac{d}{2 r}$, we have
$$
\left[H^r( \mathcal{X}  )\right]^\alpha=H^{r \alpha}( \mathcal{X}  ) \hookrightarrow C^{0, \theta}( \mathcal{X}  ) \hookrightarrow L^{\infty}( \mathcal{X}  ).
$$

Therefore $\alpha_0 = \frac{1}{\beta}$ for the embedding index of a Sobolev RKHS.

\subsubsection{dot-product kernel}

 Let $k$ be a dot-product kernel on $X = S ^d$, the unit sphere in $R ^{d+1}$, and $\mu=\sigma$ be the uniform measure on $S ^d$. Then, it is well-known that $k$ can be decomposed as
\begin{equation}
    k(x, y)=\sum_{n=0}^{\infty} \mu_n \sum_{l=1}^{a_n} Y_{n, l}(x) Y_{n, l}(y),
\end{equation}

where $\left\{Y_{n, l}\right\}$ is a set of orthonormal basis of $L^2\left( S ^d, \sigma\right)$ called the spherical harmonics. $a_n$ is multiplicity and satisfies

$$
a_n:= \binom{n+d}{n} - \binom{n-2+d}{n-2}
$$

We let $\mu_n \asymp n^{-d\beta}$ for some $\beta >1$ and also we have $a_n \asymp n^{d-1}$ and $\sum_{i=1}^{n}a_i \asymp n^{d}$, we prove embedding index $\alpha_0 = \frac{1}{\beta}$.
\begin{proof}
    By the Theorem 9 in \cite{fischer2020sobolev}, we only need to prove that for any $\alpha > \frac{1}{\beta}$, $\sum_{n=0}^{\infty} \mu_n^\alpha \sum_{l=1}^{a_n} Y_{n, l}(x)^2<\infty$.
    $$
    \begin{aligned}
\sum_{n=0}^{\infty} \mu_n^\alpha \sum_{l=1}^{a_n} Y_{n, l}(x)^2 &\le \sum_{n=0}^{\infty} \mu_n^\alpha a_n \\
&\le  \sum_{n=0}^{\infty} C n^{d-1} n^{-\alpha d \beta}\\
& = C \sum_{n=0}^{\infty} n^{-1-d(\alpha \beta-1)}<\infty
\end{aligned}
    $$
\end{proof}

\section{Detailed discussion for $f$ with complex structures}

In this section, we illustrate some $f^*$ with the complex structure and analyze the feasibility of our theories and method on these cases:

\paragraph{Mixed smoothness:} Suppose $f^*= \sum_{j=1}^{\infty}f_je_j(x)$ and $\{f_j\}_{j\in N}$ has different different decay rates. A simple example is the two-smoothness case, where $\{f_j\}_{j\in N}$ has two different decay rates, where $f_j^2 = j^{-s_1\beta-1}$ when $j$ is odd and $f_j^2 = j^{-s_2\beta-1}$ when $j$ is even ($s_1>s_2$). For the kernel $K$ with EDR $\beta$, by the definition of the interpolation space, we have
\begin{align}
    \|f^*\|_{[\mathcal{H}]^s} =\sum_{j=1}^{\infty} \frac{f_j^2}{\lambda_j^s} =\sum_{\{j: j \text{ is odd}\}} j^{s\beta -s_1\beta-1} + \sum_{\{j: j \text{ is even}\}} j^{s\beta -s_2\beta-1}
\end{align}

Thus, for $s< s_2$, we have $\|f^*\|_{[\mathcal{H}]^s}$ is bounded, meaning that $f^*_\rho\in[\mathcal{H}]^{s}$ where $s$ can be arbitrary close to $s_2$. In this case, our theory can still be applied to find the generalization ability of the kernel classifiers ($n^{-\frac{s_2\beta}{s_2\beta+1}}$). This can be also applied to multi-smoothness cases. 

However, in this case, Truncation Estimation, introduced in Section 5, is inaccurate. With a sufficient sample size, Truncation Estimation will find the smoothness $s$ between $s_1$ and $s_2$ for the two-smoothness case while we need to find $s=s_2$. In this case, the method can be improved by performing linear regression on top of $\tilde{p}_j = \sup_{k\geq j } p_j$ even though this improvement method tends to underestimate $s$. We will find new $s$ estimation methods for general situations in the near future.

\paragraph{Sobolev space of low intrinsic dimensionality}

There is a popular assumption on the real data, called manifold assumption, assuming that $f^*$ is supported on a submanifold. More specifically, for $x\in\Omega \subset \mathbb{R}^{D}$, they assume that $f^*$ belongs to the space of the low intrinsic dimensionality $d$ and $d<D$. In this case, \cite{hamm2021adaptive,ding2023random} have come up with some definitions of the low intrinsic dimension assumption:
\begin{assumption}[Low intrinsic dimension]
    There exist positive constants $c_1$ and $d\leq D$ such that for all $\delta\in (0,1)$, we have
    \begin{align}
        \mathcal{N}_{l_\infty^D}(\delta, \Omega)\leq c\delta^{-d}
    \end{align}
    where $l_\infty^D$ is  the $\mathbb{R}^D$ space equipped with $l_\infty$ norm and $\mathcal{N}_{l_\infty^D}(\delta, \Omega)$ is the covering number.
\end{assumption}

On the Sobolev space with smoothness $r$ With the assumption of low intrinsic dimension, \cite{hamm2021adaptive,ding2023random} improved their results from $n^{-\frac{2r}{2r+D}}$ to $n^{-\frac{2r}{2r+d}}$ (regression problems). Though our theories can not solve this case, we believe that the technology can be applied to classification problems and this will be our future work.

\paragraph{Well-separated data}
The well-separated assumption is another popular assumption on the real data (like MNIST, CIFAR-10, and so on) since the testing accuracy of some neural network models is near $100\%$. The well-separated assumption, in our settings, means that $f^*(x)\in\{1,-1\}$, violating the continuity of $f^*$. However, we can use a continuous function to approximate such a discontinuous function. For example, $x\in[0,1]$ and $f^*(x)=1$ if $x>1/2$ and $f^*(x)=-1$ if $x\leq 1/2$ (two regions). Then we can use an infinitely differentiable function ($s=\infty$) to approximate $f^*$ and thus the estimator finds out that the function is arbitrarily smooth. This idea can be extended to the cases with finite regions.

However, for the real data, normally with a super large dimension, the number of regions may depend on the dimension. In this situation, our theories need more effort to explain the generalization ability of the kernel classifier (like extending our theories to the high dimensional settings).

\end{document}